\documentclass[12pt,a4]{article}
\usepackage{amsmath}
\usepackage{amsfonts}
\usepackage{bm}
\usepackage{cases}
\usepackage{graphicx}
\usepackage{stmaryrd}
\usepackage{amssymb}
\usepackage{ascmac}
\usepackage{comment}
\usepackage{latexsym}
\usepackage{ascmac}
\usepackage{mathrsfs}
\usepackage{bm}			
\usepackage{comment}
\usepackage{latexsym}

\newtheorem{thm}{Theorem}[section]
\newtheorem{lem}[thm]{Lemma}
\newtheorem{prop}[thm]{Proposition}

\newtheorem{rem}[thm]{Remark}

\numberwithin{equation}{section}

\makeatletter
\newcommand{\opnorm}{\@ifstar\@opnorms\@opnorm}
\newcommand{\@opnorms}[1]{%
  \left|\mkern-1.5mu\left|\mkern-1.5mu\left|
   #1
  \right|\mkern-1.5mu\right|\mkern-1.5mu\right|
}
\newcommand{\@opnorm}[2][]{%
  \mathopen{#1|\mkern-1.5mu#1|\mkern-1.5mu#1|}
  #2
  \mathclose{#1|\mkern-1.5mu#1|\mkern-1.5mu#1|}
}
\makeatother

\newcommand{\del}{\partial}
\newcommand{\ov}{\overline}
\newcommand{\delx}{\partial_{x}}
\newcommand{\delt}{\partial_{t}}
\newcommand{\ts}{\textstyle}
\newcommand{\real}{{\rm Re}\,}

\renewcommand{\div}{\mbox{\rm div}\,}
\newcommand{\lang}{{\langle}}
\newcommand{\rang}{{\rangle}}

\newcommand{\re}{{\rm Re}} 
\newcommand{\im}{{\rm Im}}
\newcommand{\trans}{{}^\top}
\renewcommand{\Pr}{{\rm Pr}}

\newcommand{\Ker}{{\rm Ker}\,}

\newcommand{\bv}{\boldsymbol{v}}
\newcommand{\bw}{\boldsymbol{w}}
\newcommand{\bu}{\boldsymbol{u}}

\renewcommand{\b}{\boldsymbol}

\newcommand{\ep}{\varepsilon}

\renewcommand{\trans}{{}^{\top}}
\renewcommand{\lang}{{\langle}}
\renewcommand{\rang}{{\rangle}}

\newcommand{\ds}{\displaystyle}

\begin{document}

\title{\large \bf Singular limit in Hopf bifurcation for doubly diffusive convection equations I: linearized analysis at criticality 
%Spectrum of the linearized problem for a singularly perturbed system of thermal convection equations
}
\author{\normalsize Chun-Hsiung Hsia${}^1$, Yoshiyuki Kagei${}^2$, Takaaki Nishida${}^3$ \\[1ex] 
\normalsize and Yuka Teramoto${}^4$}
\date{}

\footnotetext[1]{Department of Mathematics, National Taiwan University, Taipei  10617, Taiwan}
\footnotetext[2]{Department of Mathematics, Tokyo Institute of Technology, Tokyo   152-8551, Japan}
\footnotetext[3]{Department of Advanced Mathematical Sciences, Kyoto University, Kyoto  606-8317, Japan}
\footnotetext[4]{Research Institute for Interdisciplinary Science, Okayama University, Okayama 700-8530, Japan}

\maketitle

\begin{abstract}
\noindent
A singularly perturbed system for doubly diffusive convection equations, called the artificial compressible system, is considered on a two-dimensional infinite layer for a parameters range where the Hopf bifurcation occurs in the corresponding incompressible system. The spectrum of the linearized operator in a time periodic function space is investigated in detail near the bifurcation point when the singular perturbation parameter is small. The results of this paper are the basis of the  study of the nonlinear Hopf bifurcation problem and the singular limit of the time periodic bifurcating solutions. 
\end{abstract}

\section{Introduction}\label{Introduction}
This paper is concerned with a singular limit problem for time periodic bifurcating solutions of the following system of equations:   
\begin{align}
\ep^2\delt \phi + \div\bw &= 0, \label{1.1} \\ 
\delt\bw-\Pr\Delta\bw + \Pr\nabla \phi -\Pr\mathcal{R}_1\theta\boldsymbol{e}_2 + \Pr\mathcal{R}_2\psi\boldsymbol{e}_2 + \bw\cdot\nabla\bw &= \boldsymbol{0}, \label{1.2} \\
\delt\theta - \Delta\theta - \mathcal{R}_1\bw\cdot\boldsymbol{e}_2 + \bw\cdot\nabla\theta &= 0, \label{1.3}\\
\delt \psi - d \Delta \psi - \mathcal{R}_2\bw\cdot\boldsymbol{e}_2 + \bw\cdot\nabla \psi &= 0. \label{1.4}
\end{align}
The system \eqref{1.1}--\eqref{1.4} with $\ep > 0$ is called the artificial compressible system which is a singularly perturbed system of a viscous incompressible system that describes a thermal convection phenomenon in the presence of the diffusion of salinity concentration. Here we consider the two-dimensional problem. In \eqref{1.1}--\eqref{1.4}, the unknowns $\phi = \phi(x, t)$, $\bw=\trans(w^1(x,t), w^2(x,t))$, $\theta = \theta(x, t)$ and $\psi = \psi(x, t)$ denote the deviations of the pressure, velocity field, temperature and solute concentration, respectively, at position $x \in \mathbb{R}^{2}$ and time $ t \in \mathbb{R}$, from their values of the motionless state in a thermal convection. Here and in what follows, the superscript $\trans \, \cdot $ stands for the transposition. 

The system \eqref{1.1}--\eqref{1.4} contains non-dimensional positive parameters $\mathcal{R}_{1}$, $\mathcal{R}_{2}$, $\Pr$, $d$ and $\ep$; $\mathcal{R}_{1}^{2}$ and $\mathcal{R}_{2}^{2}$ are called the thermal and salinity Rayleigh numbers, respectively; $\Pr$, $d$ and $\ep$ are called the Prandtl, Lewis and the artificial Mach numbers, respectively. The terms with $\mathcal{R}_{1}$ form a symmetric operator which may cause instabilities against the dissipative terms with Laplacian when $\mathcal{R}_{1}$ increases, while the terms with $\mathcal{R}_{2}$ form a skew-symmetric operators which may cause oscillatory behavior. The term with $\ep$ is a singular perturbation term; In the singular limit $\ep \to 0$ one obtains the incompressible system.  

We consider \eqref{1.1}--\eqref{1.4} on the two dimensional infinite layer 
\[
\{x =(x_{1}, x_{2}) \in \mathbb{R}^{2}; x_{1} \in \mathbb{R}, \, 0 < x_{2} < 1\}
\] 
under the following slip type boundary condition on the boundary $\{x_2=0,1\}$:  
\begin{equation}\label{1.5}
\mbox{$\ds \frac{\del w^1}{\del x_2} = w^2 = \theta = \psi = 0$ on $\{x_2 = 0,1\}$}.
\end{equation}
We impose the periodic boundary condition in $x_{1}$ on $\phi$, $\bw$, $\theta$ and $ \psi$ with period, say, $\frac{2\pi}{\alpha}$. We thus consider \eqref{1.1}--\eqref{1.4} in the domain 
\[
\Omega = \mathbb{T}_{\frac{2\pi}{\alpha} } \times (0, 1)
\]
under the boundary condition \eqref{1.5}. Here and in what follows, $\mathbb{T}_{a}$ denotes $\mathbb{T}_a=\mathbb{R}/a\mathbb{Z}$. 

As is mentioned above, when $\ep = 0$ in \eqref{1.1}, one obtains $\div \bw = 0$, so the system is reduced to the viscous incompressible system. The limiting process $\ep \to 0$ is mathematically formulated as a singular limit from a hyperbolic-parabolic system to a parabolic system. It is thus of interest to consider the question whether the artificial compressible system approximates well the incompressible system in the limit $\ep \to 0$. 

The aim of this paper is to investigate whether the artificial compressible system \eqref{1.1}--\eqref{1.4} gives a good approximation of the incompressible system $\eqref{1.1}$--\eqref{1.4} with $\ep = 0$ as $\ep \to 0$ when a time periodic bifurcation occurs in the incompressible system $\eqref{1.1}$ --\eqref{1.4} with $\ep = 0$.  

We briefly review bifurcation results for the incompressible convection system. Let $\ep = 0$ in \eqref{1.1}--\eqref{1.4}. We denote by $u = \trans(\phi, \bw, \theta, \psi)$ the perturbation of the motionless state. If $\mathcal{R}_{2} =0$, the equation \eqref{1.4} for the evolution of the solute concentration $\psi$ becomes a simple convection-diffusion equation and essentially decouples with $\eqref{1.1}$--\eqref{1.3} with $\ep = 0$. So, in this case, the convection phenomenon is essentially governed by $\eqref{1.1}$--\eqref{1.3} with $\ep = 0$ for the unknowns $\trans(\phi, \bw, \theta)$, which is called the Oberbeck-Boussinesq system. As to the Oberbeck-Boussinesq system, it is well known (\cite{Busse, Iudovich, Joseph, Rabinowitz}) that if $\mathcal{R}_{1}$ is sufficiently small, then the motionless state $u = 0$ is asymptotically stable. When $\mathcal{R}_{1}$ increases, the motionless state $u = 0$ becomes unstable beyond a certain critical value of $\mathcal{R}_{1}$ and stationary convective patterns bifurcate from $u = 0$ for some range of $\alpha$. In the case of the Oberbeck-Boussinesq system, the onset of convection is caused only by stationary bifurcation due to a variational structure.  
On the other hand, when $\mathcal{R}_{2} > 0$, the presence of the skew-symmetric effect enriches the dynamics near the onset of convection. Indeed, it was shown in \cite{Bona-Hsia-Ma-Wang} that for some range of $\mathcal{R}_{2}$, $d$ and $\alpha$, a Hopf bifurcation occurs, namely, time periodic convective patterns bifurcate from $u = 0$ beyond a critical value of $\mathcal{R}_{1}$.  For more details about the convection problem with salinity, see \cite{Bona-Hsia-Ma-Wang, Hsia-Ma-Wang} and references therein. 

The artificial compressible system for the incompressible Navier-Stokes equations 
\begin{align} 
\ep^{2}\delt p + \div \bv & = 0, \label{1.6} \\
\delt \bv - \nu \Delta \bv + \nabla p + \bv \cdot \nabla \bv & = \b{0} \label{1.7}
\end{align} 
was proposed by Chorin (\cite{Chorin1, Chorin2, Chorin3}) and Temam (\cite{Temam1,Temam2}) to avoid the difficulties in numerical computations caused by the constraint $\eqref{1.6}$ with $\ep = 0$, i.e., $\div \bv = 0$. By using the artificial compressible system for the Oberbeck-Boussinesq system \eqref{1.1}--\eqref{1.3} with $\mathcal{R}_{2} = 0$, Chorin computed stationary periodic convective patterns near the onset of convection. 
The convergence of solutions as $\ep \to 0$ was discussed by Temam in \cite{Temam1, Temam2, Temam3} for the artificial compressible system with the additional stabilizing nonlinear term $+\frac{1}{2}(\div \bv)\bv$ on the left-hand side of \eqref{1.7}: 
\begin{align}\label{1.8}
\delt \bv - \nu \Delta \bv + \nabla p + \bv \cdot \nabla \bv +\frac{1}{2}(\div \bv)\bv & = \b{0}.  
\end{align}
It was proved in \cite{Temam1, Temam2, Temam3} that the solutions $\{p^{\ep}, \bv^{\ep}\}$ of  the initial boundary value problem for the artificial compressible system \eqref{1.6}, \eqref{1.8} on a two-dimensional bounded domain $D$ converge to the one for the incompressible system as $\ep \to 0$ in such a way that 
$\bv_{\ep} \to \bv$ strongly in $L^2(0,T;H^{1}(D)^2)$ and $\nabla p_{\ep} \to \nabla p$ strongly in $H^{-1}(\Omega \times (0,T))$ for all $T>0$, where $\trans(p,\bv)$ is the solution of the initial boundary value problem for the incompressible Navier-Stokes equations \eqref{1.6}--\eqref{1.7} with $\ep = 0$. Convergence result was also established in the framework of weak solutions on three-dimensional bounded domains. In this direction, Donatelli \cite{Donatelli1, Donatelli2} and Donatelli and Marcati \cite{Donatelli-Marcati1, Donatelli-Marcati2} proved convergence results which reflect a dispersive aspect of the system in the case of unbounded domains by using the wave equation structure of the pressure and the dispersive estimates. 

In \cite{Kagei-Nishida, Kagei-Nishida-Teramoto, Teramoto}, it was investigated whether the artificial compressible systems \eqref{1.6}--\eqref{1.7} and \eqref{1.1}--\eqref{1.3} with $\mathcal{R}_{2} = 0$ give good approximations of the corresponding incompressible systems as $\ep \to 0$ from the view point of the stability of stationary solutions. For the stability questions, one needs to investigate the spectrum of the linearized operators around a stationary solution. It was shown in \cite{Teramoto} that if the stationary bifurcation occurs in the incompressible system, then a part of the spectrum near the origin of the linearized operator around the bifurcating  stationary solution is approximated by a part of the spectra near the origin  of the corresponding linearized operators for the artificial compressible system as $\ep \to 0$. Furthermore, if the basic flow satisfies the energy-type stability criterion given in \cite{Kagei-Nishida-Teramoto}, the stability of bifurcating stationary solutions of both the incompressible and artificial compressible systems coincide with each other. The result is applicable to the bifurcating convective patterns near the onset of convection for \eqref{1.1}--\eqref{1.4}, and hence, one can conclude that the artificial compressible system \eqref{1.1}--\eqref{1.4} gives a good approximation of the incompressible system in the limit $\ep \to 0$ near the onset of convection caused by stationary bifurcation in the incompressible system \eqref{1.1}--\eqref{1.4} with $\ep = 0$. 

We thus raise the question how about the case of Hopf bifurcation. The main difference to the case of stationary bifurcation is as follows. The sets of stationary solutions of the artificial compressible system and the incompressible system are the same, so one does not need to show the occurrence of stationary bifurcation. In contrast to the stationary case, in the case of the Hopf bifurcation, one first needs to show the occurrence of a Hopf bifurcation for  \eqref{1.1}--\eqref{1.4} with uniform estimates in $0 < \ep \ll 1$, which requires us a more detailed analysis of the spectral properties of the linearized semigroup around the basic flow. This is of course due to that the singular limit under question is of type of vanishing time derivatives. In this paper, we study the linearized problem around $u = 0$ of \eqref{1.1}--\eqref{1.4} with $0 < \ep \ll 1$ near the bifurcation point of time periodic bifurcation in the incompressible problem \eqref{1.1}--\eqref{1.4} with $\ep = 0$. Based on the analysis of this paper, the nonlinear bifurcation problem is studied in the paper \cite{Hsia-Kagei-Nishida-Teramoto2}, where it is proved that the Hopf bifurcation occurs in the artificial compressible system \eqref{1.1}--\eqref{1.4} for $0 < \ep \ll 1$, together with the convergence of the time periodic bifurcating branch as $\ep \to 0$. 

To explain the results of this paper more precisely, we write the time periodic problem for \eqref{1.1}--\eqref{1.4} in the form  
\begin{equation}\label{1.9}
\delt u + L^{\ep}_{\mathcal{R}_{1} } u + N(u) = 0,
\end{equation}
where $N(u)$ denotes the nonlinearity; and $L^{\ep}_{\mathcal{R}_{1} }$ denotes the linearized operator around $u = 0$: 
\[
L^{\ep}_{\mathcal{R}_1} =
\begin{pmatrix}
0 & \frac{1}{\ep^2}\div & 0 & 0\\
\Pr\nabla & -\Pr\Delta & -\Pr\,\mathcal{R}_{1}\boldsymbol{e}_2 & \Pr\mathcal{R}_2\boldsymbol{e}_2\\
0 & -\mathcal{R}_{1}\trans\boldsymbol{e}_2  & -\Delta & 0\\
0 &-\mathcal{R}_2\trans\boldsymbol{e}_2 & 0 & - d \Delta
\end{pmatrix}.
\]
Here $\b{e}_{2} = \trans(0, 1)$. We take $\mathcal{R}_1$ as a bifurcation parameter. The precise functional setting is given in section \ref{incompressible}. 

We denote by $\mathbb{L}_{\mathcal{R}_{1}}$ the corresponding linearized operator for the incompressible problem \eqref{1.1}--\eqref{1.4} with $\ep = 0$. It was shown in \cite{Bona-Hsia-Ma-Wang} that, for $(\mathcal{R}_{2}, d)$ in a certain range, there exists a critical value $\mathcal{R}_{1,c}$ such that when $\mathcal{R}_{1}$ increases, a pair of complex conjugate simple eigenvalues $\{\lambda_{+}(\mathcal{R}_{1}), \lambda_{-}(\mathcal{R}_{1})\}$ with $\ov{\lambda_{-}(\mathcal{R}_{1})} = \lambda_{+}(\mathcal{R}_{1})$ of $- \mathbb{L}_{\mathcal{R}_{1}}$ crosses the imaginary axis at $\pm ia$ for $\mathcal{R}_{1} = \mathcal{R}_{1,c}$ and the remaining part of the spectrum of $ - \mathbb{L}_{\mathcal{R}_{1, c}}$ remains in the left-half plane strictly away from the imaginary axis. 

We shall show that if $\ep$ is sufficiently small, then there exists a critical value $\mathcal{R}^{\ep}_{1,c} = \mathcal{R}_{1,c} + O(\ep^{2})$ such that when $\mathcal{R}_{1}$ increases, a pair of complex conjugate simple eigenvalues $\{\lambda^{\ep}_{+}(\mathcal{R}_{1}), \lambda^{\ep}_{-}(\mathcal{R}_{1})\}$ with $\ov{\lambda^{\ep}_{-}(\mathcal{R}_{1})} = \lambda^{\ep}_{+}(\mathcal{R}_{1}) = \lambda_{+}(\mathcal{R}_{1}) + O(\ep^{2})$ of $- L^{\ep}_{\mathcal{R}_{1}}$ crosses the imaginary axis at $\pm ia^{\ep} = \pm i(a + O(\ep^{2}) )$ for $\mathcal{R}_{1} = \mathcal{R}^{\ep}_{1,c}$ and the remaining part of the spectrum of $ - L^{\ep}_{\mathcal{R}^{\ep}_{1, c}}$ remains in the left-half plain strictly away from the imaginary axis uniformly in small $\ep$. This result is proved by combining the perturbation arguments in  \cite{Kagei-Nishida, Kagei-Teramoto}. 

We shall then solve the linearized time periodic problem 
\[
\delt u + L^{\ep}_{\mathcal{R}_{1} } u = F 
\]
by the Lyapunov-Schmidt method for $\mathcal{R}_{1}$ near the critical value $\mathcal{R}^{\ep}_{1,c}$. Here $F$ is a given time periodic function. 
For this purpose, we investigate the spectrum of $\tilde{B}^{\ep} = \delt + L^{\ep}_{\mathcal{R}^{\ep}_{1,c} }$ on a space $\mathcal{X}_{a}$ of time periodic functions. We first study the spectral properties of the linearized semigroup $e^{ - t L^{\ep}_{\mathcal{R}^{\ep}_{1,c} } }$. A key in the analysis is to establish an exponential decay estimate for the semigroup $e^{- t L^{\ep}_{\mathcal{R}^{\ep}_{1,c} } }$ in the complementary subspace to the eigenspaces for the eigenvalues $\{\lambda^{\ep}_{+}(\mathcal{R}^{\ep}_{1,c}), \lambda^{\ep}_{-}(\mathcal{R}^{\ep}_{1,c})\}$ uniformly in small $\ep$. To this end, we employ the energy method with $\ep$ weights and make use of an oscillatory aspect of the semigroup. 
 
Based on the spectral properties of $e^{ - t L^{\ep}_{\mathcal{R}^{\ep}_{1,c} } }$, we deduce that 
\[
\rho(-\tilde{B}^{\ep})\supset \Sigma_1 \setminus (\cup_{k \in \mathbb{Z}} \{ika^{\ep} \}),
\]
where $\Sigma_1=\{\lambda\in\mathbb{C};\re\lambda>- \kappa_1 \}$ with some positive constant $\kappa_1$ independent of small $\ep$; each $ika^{\ep} $ is a semisimple eigenvalue of $\tilde{B}^{\ep}$ with two-dimensional eigenspace. In particular, we have that $0$ is a semisimple eigenvalue of $\tilde{B}^{\ep}$ and $\mathcal{X}_{a} = \Ker (\tilde{B}^{\ep}) \oplus R(\tilde{B}^{\ep})$. Furthermore, we establish a uniform estimate of the inverse $(\tilde{B}^{\ep}|_{R(\tilde{B}^{\ep})})^{-1}$ with respect to $0 < \ep \ll 1$. 

From these spectral properties of $\tilde{B}^{\ep}$, one could expect that a Hopf bifurcation would occur in \eqref{1.1}--\eqref{1.4} for $\mathcal{R}_{1}$ near $\mathcal{R}^{\ep}_{1,c}$. Indeed, one can show the existence of a nontrivial time periodic solution branch by a variance of the standard bifurcation theory. However, we need to establish uniform estimates in $\ep$ for bifurcating branch to consider its convergence to the bifurcating branch of the incompressible system. For this purpose, we investigate the spectral properties of the operator $(1 + \omega)\delt + L^{\ep}_{\mathcal{R}^{\ep}_{1,c} }$ for small $\ep$ and $\omega$ in a time periodic function space, which is done by using a uniform exponential decay estimate in $\ep$ and $\omega$ for $e^{-\frac{ t}{ 1 + \omega} L^{\ep}_{\mathcal{R}^{\ep}_{1,c} } }$ in the complementary subspace to the eigenspaces for the eigenvalues $\{\lambda^{\ep}_{+}(\mathcal{R}^{\ep}_{1,c}), \lambda^{\ep}_{-}(\mathcal{R}^{\ep}_{1,c})\}$. In contrast to Hopf bifurcation problem for the incompressible equations, $\omega\delt$ cannot be regarded as a simple perturbation of the operator $\delt + L^{\ep}_{\mathcal{R}^{\ep}_{1,c} }$ even if $\omega$ is sufficiently small; $\omega\delt$ could cause ``$\ep$-loss'' if it would be regarded as a perturbation of $\delt + L^{\ep}_{\mathcal{R}^{\ep}_{1,c} }$. 

This paper is organized as follows. In section 2 we introduce notation used in this paper. In section 3 we state the results on the spectrum of the linearized operator for the incompressible problem obtained in \cite{Bona-Hsia-Ma-Wang}. 
In section 4 we state the main results of this paper, i.e., the results on the spectra of $- L^{\ep}_{\mathcal{R}^{\ep}_{1,c} }$ and $(1 + \omega)\delt + L^{\ep}_{\mathcal{R}^{\ep}_{1,c} }$. 
Section 5 is devoted to investigating the spectral properties of the operator $(1 + \omega)\delt + L^{\ep}_{\mathcal{R}^{\ep}_{1,c} }$ for $0 < \ep \ll 1
$ and $|\omega| \ll 1$; and we give proofs of the main results. 

\section{Preliminaries}

In this section we introduce notation used in this paper. 
Let $\Omega=\mathbb{T}_{\frac{2\pi}{\alpha}}\times(0,1)$. The usual $L^p$ space on $\Omega$ is denoted by $L^p(\Omega)$ with norm $\|\cdot\|_p$. The space of all $L^p$ vector fields on $\Omega$ is denoted by $\b L^p(\Omega)$ and its norm is also denoted by $\|\cdot\|_p$. 
The inner product of $L^2(\Omega)$ is denoted by $(\cdot, \cdot)_{L^2}$, and the inner product of $\b L^2(\Omega)$ is also denoted by $(\cdot,\cdot)_{L^2}$. 

For $\bu_j=\trans(\bw_j,\theta_j,\psi_j)\in \b L^2(\Omega)\times L^2(\Omega)\times L^2(\Omega)$ $(j=1,2)$, we define the inner product $(\bu_1,\bu_2)$ by 
\[
(\bu_1,\bu_2)=\Pr^{-1}(\bw_1,\bw_2)_{L^2}+(\theta_1,\theta_2)_{L^2}+(\psi_1,\psi_2)_{L^2} 
\]
and the norm $\|\bu\|_2$ of $\bu=\trans(\bw,\theta,\psi)\in \b L^2(\Omega)\times L^2(\Omega)\times L^2(\Omega)$ by 
\[
\|\bu\|_2=\sqrt{(\bu,\bu)}=\sqrt{\Pr^{-1}\|\bw\|_2^2+\|\theta\|_2^2+\|\psi\|_2^2}. 
\]
For $\bu_j=\trans(\bw_j,\theta_j,\psi_j)\in L^2(\mathbb{T}_{\frac{2\pi}{a}};\b L^2(\Omega)\times L^2(\Omega)\times L^2(\Omega))$ $(j=1,2)$, we define the inner product $\lang \bu_1,\bu_2\rang$ by 
\[
\langle\bu_1,\bu_2\rangle = \frac{a}{2\pi}\int_0^{\frac{2\pi}{a}}(\bu_1(t),\bu_2(t))\,dt.  
\]

Let $\ep$ be a given positive number. 
For $u_j=\trans(\phi_j,\bu_j)\in L^2(\Omega)\times \b L^2(\Omega)\times L^2(\Omega)\times L^2(\Omega)$ $(j=1,2)$, we define the inner product $(u_1,u_2)_\ep$ by
\[
(u_1,u_2)_\ep = \ep^2(\phi_1,\phi_2)_{L^2} + (\bu_1,\bu_2)
\]
and the norm $|||u|||_\ep$ of $u=\trans(\phi,\bu)\in  L^2(\Omega)\times \b L^2(\Omega)\times L^2(\Omega)\times L^2(\Omega)$ by 
\[
|||u|||_\ep=\sqrt{(u,u)_\ep}=\sqrt{\ep^2\|\phi\|_2^2+\|\bu\|_2^2},
\]
and, likewise, for $\bu_j=\trans(\bw_j,\theta_j,\psi_j)\in L^2(\mathbb{T}_{\frac{2\pi}{a}};L^2(\Omega)\times \b L^2(\Omega)\times L^2(\Omega)\times L^2(\Omega))$ $(j=1,2)$, the inner product $\langle u_1,u_2\rangle_\ep$ is defined by  
\[
\langle u_1,u_2\rangle_\ep = \frac{a}{2\pi}\int_0^{\frac{2\pi}{a}}(u_1(t),u_2(t))_\ep\,dt. 
\] 

The $k$th order $L^2$ Sobolev spaces are written as $H^k(\Omega)$ for scalar functions and $\b H^k(\Omega)$ for vector fields. 
We set 
\[
L_*^2(\Omega)=\{\phi \in L^2(\Omega); \int_{\Omega}\phi(x) \,dx=0\}, 
\ \ \ 
H_*^k(\Omega)=H^k(\Omega)\cap L_*^2(\Omega). 
\]

We next introduce function spaces with symmetries. We set 
\begin{align*}
L_{sym}^2(\Omega) & = \{\theta\in L^2(\Omega); \mbox{ \rm $\theta$ is even in $x_1$ } \},
\\
L_{*,sym}^2(\Omega)& =L_{sym}^2(\Omega) \cap L_*^2(\Omega), 
\\
H^k_{*,sym}(\Omega)& =H^k(\Omega)\cap L_{*,sym}^2(\Omega) 
\end{align*}
and 
\[
\b L^2_{sym}(\Omega)=\{\bv = \trans(v^1,v^2)\in \b L^2(\Omega); \mbox{ \rm $v^1$ is odd in $x_1$, $v^2$ is even in $x_1$} \}. 
\]
It is known (\cite{Galdi, Sohr, Temam3}) that $\b L^2_{sym}(\Omega)$ admits the Helmholtz decomposition: 
\[
\mbox{$\b L^2_{sym}(\Omega)=\b L^2_{\sigma,sym}(\Omega)\oplus \b G_{sym}^2(\Omega)$ (orthogonal decomposition)}, 
\]
where   
\[
\b L^2_{\sigma,sym}(\Omega)
=\{\bw \in\b L_{sym}^2(\Omega);\,\div\bw=0\mbox{ in }\Omega,\,\bw\cdot\b n|_{\del\Omega}=0\}
\]
and 
\[
\b G_{sym}^2(\Omega)=\{\nabla \phi; \phi \in H_{*,sym}^1(\Omega)\}.
\]
The orthogonal projection on $\b L^2_{\sigma,sym}(\Omega)$ is denoted by $\mathbb{P}_\sigma$. 

We also introduce the following function spaces: 
\begin{align*}
\b X &= \b L^2_{sym}(\Omega) \times L_{sym}^2(\Omega)\times L_{sym}^2(\Omega), 
\\
\b Y & =[\b H_{b}^2(\Omega)\times (H^2\cap H^1_0)(\Omega)\times (H^2\cap H^1_0)(\Omega))]\cap \b X, 
\\ 
\mathbb{X}_\sigma &= \b L^2_{\sigma,sym}(\Omega) \times L_{sym}^2(\Omega)\times L_{sym}^2(\Omega), 
\\ 
\mathbb{Y}_\sigma & =[\b H_{b}^2(\Omega)\times (H^2\cap H^1_0)(\Omega)\times (H^2\cap H^1_0)(\Omega))]\cap \mathbb{X}_\sigma, 
\\
X & =H_{*,sym}^1(\Omega)\times \b L^2_{sym}(\Omega) \times L_{sym}^2(\Omega)\times L_{sym}^2(\Omega),  
\\
Y & =[H_{*,sym}^1(\Omega)\times \b H_{b}^2(\Omega)\times (H^2\cap H^1_0)(\Omega)\times (H^2\cap H^1_0)(\Omega))]\cap X.
\end{align*}
Here 
\[
\b H_b^2(\Omega) = \left\{\bw = \trans(w^1,w^2)\in \b H^2(\Omega);\,\frac{\del w^1}{\del x_2} = w^2 = 0\mbox{ on }\{x_2 = 0,1\}\right\}
\]
and 
\[
H_0^1(\Omega) = \{\theta\in H^1(\Omega);\,\theta|_{x_2 = 0,1} = 0\}.
\]

We also introduce the space $\b X^1$ defined by 
\[
\b X^1=[\b H_b^1(\Omega)\times H_0^1(\Omega)\times H_0^1(\Omega)]\cap \b X, 
\]
where 
\[
\b H_b^1(\Omega) = \left\{\bw = \trans(w^1,w^2)\in \b H^1(\Omega);\, w^2|_{x_2 = 0, 1} = 0\right\}.
\]
For $\bu=\trans(\bw,\theta,\psi) \in \b X^1$, we set 
\[
\|\bu\|_{\b X^1}=\sqrt{\Pr^{-1}\|\nabla\bw\|_2^2+\|\nabla\theta\|_2^2+\|\nabla\psi\|_2^2}
\]
We note that the Poincar\'e inequality $\|\bw\|_2\le C\|\nabla\bw\|_2$ holds for $\bw=\trans(w^1,w^2) \in \b H_b^1(\Omega)\cap \b L^2_{sym}(\Omega)$ since $\int_\Omega w^1\,dx=0$ because of the oddness of $w^1$ in $x_1$. Therefore, $\|\cdot\|_{\b X^1}$ defines a norm on $\b X^1$.  

If $\bu=\trans(\bw,\theta,\psi) \in \b X$ with $\bw=\trans(w^1,w^2)$, then $\bu$ is expanded as Fourier series: 
\begin{align*}
w^1 & =\sum_{j\ge1, k\ge0}w_{jk}^1\sin{\alpha j x_1}\cos{k\pi x_2}, 
\\
w^2 & =\sum_{j\ge0, k\ge1}w_{jk}^2\cos{\alpha j x_1}\sin{k\pi x_2}, 
\\
\theta & =\sum_{j\ge0, k\ge1}\theta_{jk}\cos{\alpha j x_1}\sin{k\pi x_2}, 
\\
\psi & =\sum_{j\ge0, k\ge1}\psi_{jk}\cos{\alpha j x_1}\sin{k\pi x_2}.
\end{align*}
In terms of this expansion of $\bu\in \b X$, we define $\|\bu\|_{(\b X^1)^*}$ by 
\[
\|u\|_{(\b X^1)^*}=\left[\sum_{j\ge0,k\ge0}(\alpha^2 j^2+k^2\pi^2)^{-1}\left(\Pr^{-1}\left\{(w_{jk}^1)^2+(w_{jk}^2)^2\right\}+\theta_{jk}^2+\psi_{jk}^2 \right) \right]^{\frac{1}{2}},
\]
where $w_{0k}^1=w_{j0}^2=\theta_{j0}=\psi_{j0}=0$ for $k,j\ge0$. 
It then follows that 
\[
\| \bu \|_{(\b X^1)^*} \le C \| \bu \|_2
\]
and 
\[
(\bu_1,\bu_2) \le \|\bu_1\|_{\b X^1}\|\bu_2\|_{(\b X^1)^*}.
\]
We also define $X^1$ by 
\[
X^1 = H_{*,sym}^1(\Omega) \times \b X^1.
\] 

For given $t_{1} < t_{2}$ we set  
\begin{align*}
\mathcal{X} (t_{1}, t_{2}) & =L^2(t_{1}, t_{2}; X), 
\\
\mathcal{Y} (t_{1}, t_{2})& =L^2(t_{1}, t_{2},T;Y)\cap H^1(t_{1}, t_{2},T; X), 
\end{align*}
and the norm of $u\in \mathcal{X}(t_{1}, t_{2})$ (resp. $u\in \mathcal{Y}(t_{1}, t_{2})$) is denoted by $\|u\|_{\mathcal{X}(t_{1}, t_{2})}$ (resp. $\|u\|_{\mathcal{Y}(t_{1}, t_{2})}$). 

We shall also use the norms $||| u |||_{\ep, X^1}$, $|||u |||_{\ep, \mathcal{X}(t_{1}, t_{2})}$ and $|||u |||_{\ep, \mathcal{Y}(t_{1}, t_{2})}$ of $u=\trans(\phi, \bu)$ with $\ep$ weights defined by 
\begin{align*}
||| u |||_{\ep,X^1} & = \left\{ |||u|||_\ep^2 + \ep^2 |||\delx u|||_\ep^2 \right\}^{\frac{1}{2}},
\\
||| u |||_{\ep, \mathcal{X}(t_{1}, t_{2})}  & = \left\{\int_{t_{1} }^{t_{2} } \left( \ep^2 \| \phi \|_2^2 + \|\bu\|_{(\b X^1)^*}^2 + \ep^6 \|\delx \phi \|_2^2 +\ep^2 \|\bu\|_2^2\right) \, dt \right\}^{\frac{1}{2}}, 
\\
||| u |||_{\ep, \mathcal{Y}(t_{1}, t_{2})} & = \left\{ \sup_{t_{1} \le t \le t_{2} } ||| u(t) |||_{\ep,X^1}^2 \right. \\
& \quad \quad \left. + \int_{t_{1} }^{t_{2} } \left(|||\delx u |||_\ep^2 + \ep ^2 |||\delt u |||_\ep^2 + \ep^2 \|\delx ^2\bu \|_2^2 + \ep^6 \|\delx\delt \phi\|_2^2 \right) \,dt \right\}^{\frac{1}{2}}. 
\end{align*}

We next introduce functions spaces of time periodic functions with period $\frac{2\pi}{a}$. We define the spaces $\mathcal{X}_a$ and $\mathcal{Y}_a$ by 
\begin{align*}
\mathcal{X}_a & =L^2(\mathbb{T}_{\frac{2\pi}{a}};X), 
\\
\mathcal{Y}_a & =L^2(\mathbb{T}_{\frac{2\pi}{a}};Y)\cap H^1(\mathbb{T}_{\frac{2\pi}{a}};X), 
\end{align*}
respectively. The norm of $u\in \mathcal{X}_a$ (resp. $u\in \mathcal{Y}_a$) is defined by $\|u\|_{\mathcal{X}_a} = \|u\|_{\mathcal{X}(0, \frac{2\pi}{a})}$ (resp. $\|u\|_{\mathcal{Y}_a} = \|u\|_{\mathcal{Y}(0, \frac{2\pi}{a})}$). The weighted norm of $u\in \mathcal{X}_a$ (resp. $u\in \mathcal{Y}_a$) is defined by $||| u |||_{\ep, \mathcal{X}_a} = ||| u |||_{\ep, \mathcal{X}(0, \frac{2\pi}{a})}$ (resp. $||| u |||_{\ep, \mathcal{Y}_a} = ||| u |||_{\ep, \mathcal{Y}(0, \frac{2\pi}{a})}$). 

For an operator $A$, we denote the resolvent set of $A$ by $\rho(A)$ and the spectrum of $A$ by $\sigma(A)$. The space of all bounded linear operators from $E_{1}$ to $E_{2}$ is denoted by $\mathfrak{B}(E_{1}, E_{2})$.

\section{Spectrum of the linearized operator for the incompressible system}\label{incompressible}

In this section we recall the results on the spectrum of the linearized operator at the motionless state for the incompressible system \eqref{1.1} obtained in \cite{Bona-Hsia-Ma-Wang}. 

We first derive the non-dimensional system \eqref{1.1}--\eqref{1.4} for the perturbation of the motionless state. Consider the incompressible system 
\begin{equation}\label{3.1}
\left\{
\begin{array}{rcl}
\div\bv  &=& 0, \\ 
\delt\bv-\nu\Delta\bv + \frac{1}{\rho_0}\nabla p + \bv\cdot\nabla\bv & = &-\frac{\rho \mathrm{g}}{\rho_0} \b e_2 ,\\
\delt T - d_T\Delta T + \bv\cdot\nabla T & = & 0,\\ 
\delt S - d_S\Delta S + \bv\cdot\nabla S  & = & 0
\end{array}
\right.
\end{equation}
in $\mathbb{T}_{\frac{2\pi}{\alpha}}\times(0,\ell)$. 
Here $\alpha$ and $\ell$ are given positive constants; $p=p(x,t)$, $\bv=\trans(v^1(x,t),v^2(x,t))$, $T=T(x,t)$ and $S=S(x,t)$ are the pressure, velocity field, temperature and solute concentration, respectively, at position $x=(x_1,x_2)\in \mathbb{T}_{\frac{2\pi}{\alpha}}\times(0,\ell)$ and time $t \in \mathbb{R}$; $\mathrm{g}$ is the gravity constant; $\b e_2=\trans(0,1)$ is the unit vector in $x_2$ direction; $\nu$ is the kinematic viscosity; $d_T$ is the thermal conductivity; $d_S$ is the solute diffusivity; 
$\rho_0$ is the fluid density on the boundary $\{x_2=0\}$; $\rho=\rho(x,t)$ is the fluid density which is assumed to have the form 
\[
\rho=\rho_0\left[1-a_T(T-T_0)+a_S(S-S_0)\right].
\] 
Here $\mathrm{g}$, $\nu$, $d_T$, $d_S$, $a_T$ and $a_S$ are assumed to be positive constants.  

We consider \eqref{3.1} under the following boundary conditions on the boundaries $\{x_2=0,\ell\}$: 
\begin{equation}\label{3.2}
\begin{cases}
\del_{x_2}v^1=v^2=0, \ \ T=T_0, \ \ S=S_0 \ \ \mbox{on $\{x_2=0\}$},  \\
\del_{x_2}v^1=v^2=0, \ \ T=T_1, \ \ S=S_1 \ \ \mbox{on $\{x_2=\ell\}$}.
\end{cases}
\end{equation}
Here $T_j$ and $S_j$ $(j=0,1)$ are positive constants satisfying 
\[
\mbox{$T_0>T_1$ and $S_0>S_1$}. 
\] 

One can verify that the problem \eqref{3.1}--\eqref{3.2} has the following stationary solution $u_B=\trans(p_B,\bv_B,T_B,S_B)$: 
\begin{equation}\label{3.3}
\begin{cases}
p_B=p_0-\mathrm{g}\left[x_2-\frac{a_T(T-T_0)}{2\ell}x_2^2+\frac{a_S(S-S_0)}{2\ell}x_2^2\right], \\
\bv_B=\b 0, \\
T_B=\frac{T_1-T_0}{\ell}x_2+T_0, \\ 
S_B=\frac{S_1-S_0}{\ell}x_2+S_0,
\end{cases}
\end{equation}
where $p_0$ is a positive constant. 

We are interested in the stability of the basic state $u_B=\trans(p_B,\bv_B,T_B,S_B)$, so we begin with rewriting  the problem \eqref{3.1}--\eqref{3.2} into a non-dimensional form. 
We introduce the following non-dimensional variables: 
\[
x=\ell x_*, \ t=\frac{\ell^2}{d_T}t_*, \ p=\frac{\rho_0\nu d_T}{\ell^2}p_*, \
T=(T_0-T_1)T_*, \ S=(S_0-S_1)S_*.
\]
Under this transformation, the domain $\mathbb{T}_{\frac{2\pi}{\alpha}}\times (0,\ell)$ is transformed into 
\[
\Omega=\mathbb{T}_{\frac{2\pi}{\alpha_*}}\times (0,1),
\]
where $\alpha_*=\alpha \ell$. 
Let $u_{B*}=\trans(p_{B*},\bv_{B*},T_{B*},S_{B*})$ denote the non-dimensionalized basic state and let the perturbation be denoted by 
\[
u  =\trans(\phi,\bw,\theta,\psi) =\trans(p_*-p_{B*},\bv_*-\bv_{B*}, \mathcal{R}_1(T_*-T_{B*}), \mathcal{R}_2(S_*-S_{B*}))
\]
with $\mathcal{R}_j=\sqrt{R_j}$ $(j=1,2)$, where $R_1$ and $R_2$ are the thermal and salinity Rayleigh numbers, respectively, defined below. 

After omitting $*$ of $x_*$ and $t_*$, we see that the perturbation $u=\trans(\phi,\bw,\theta,\psi)$ is governed by the following system of equations: 
\begin{equation}\label{3.4}
\left\{
\begin{array}{rcl}
\div\bw &=& 0,\\
\delt\bw-\Pr\Delta\bw + \Pr\nabla \phi -\Pr\mathcal{R}_1\theta\b{e}_2 + \Pr\mathcal{R}_2\psi\b{e}_2 +\bw\cdot\nabla\bw &=& \b{0},\\
\delt\theta - \Delta\theta - \mathcal{R}_1\bw\cdot\b e_2 + \bw\cdot\nabla\theta &=& 0,\\
\delt \psi - d \Delta \psi - \mathcal{R}_2\bw\cdot\b e_2 + \bw\cdot\nabla \psi &=& 0.
\end{array}
\right.
\end{equation}
Here $\Pr$ and $d$ are the non-dimensional numbers defined by 
\[
\mbox{$\ds \Pr=\frac{\nu}{d_T}$: the Prandtl number, $\ds d=\frac{d_S}{d_T}$: the Lewis number}, 
\]
and $\mathcal{R}_j=\sqrt{R_j}$, $j=1,2$, are the non-dimensional numbers with $R_j$ defined by 
\[
\begin{array}{l}
\mbox{$\ds R_1=\frac{a_T\mathrm{g}\ell^3(T_0-T_1)}{d_T\nu}$: the thermal Rayleigh number}, \\[1ex]
\mbox{$\ds R_2=\frac{a_S\mathrm{g}\ell^3(S_0-S_1)}{d_T\nu}$: the salinity Rayleigh number}. 
\end{array}
\]
The boundary condition on the boundary $\{x_2=0,1\}$ is given as 
\begin{equation}\label{3.5}
\mbox{$\ds \frac{\del w^1}{\del x_2} = w^2 = \theta = \psi = 0$ on $\{x_2 = 0,1\}$.}
\end{equation}
We note that the problem \eqref{3.4}--\eqref{3.5} has a trivial stationary solution $u=0$ which corresponds to the basic state $u_B$. 

As for the incompressible problem \eqref{3.4}--\eqref{3.5}, Bona, Hsia, Ma and Wang showed in \cite{Bona-Hsia-Ma-Wang} that 
there are positive numbers $\mathcal{R}_{2*}$ and $\mathcal{R}_2^*$ such that if $\mathcal{R}_2\in[\mathcal{R}_{2*},\mathcal{R}_2^*)$, $\Pr>1$ and $0<d<1$, then there exists a critical number $\mathcal{R}_{1,c}$ such that the basic state $u=0$ is stable when $\mathcal{R}_1<\mathcal{R}_{1,c}$, while $u=0$ is unstable when $\mathcal{R}_1>\mathcal{R}_{1,c}$, and time periodic solutions bifurcate from $u=0$ for $\mathcal{R}_1>\mathcal{R}_{1,c}$. 

More precisely, we introduce a bifurcation parameter $\eta = \mathcal{R}_1-\mathcal{R}_{1,c}$. 
We define the linearized operator $\mathbb{L}_{\mathcal{R}_{1,c}+\eta}$ on $\mathbb{X}_\sigma$ by 
\[
D(\mathbb{L}_{\mathcal{R}_{1,c}+\eta})=\mathbb{Y}_\sigma,  
\]
\[
\mathbb{L}_{\mathcal{R}_{1,c}+\eta} = 
\begin{pmatrix}
-\Pr\mathbb{P}_\sigma\Delta & -\Pr(\mathcal{R}_{1,c} + \eta)\mathbb{P}_\sigma\boldsymbol{e}_2 & \Pr\mathcal{R}_2\mathbb{P}_\sigma\boldsymbol{e}_2 \\
-(\mathcal{R}_{1,c} + \eta)\trans\boldsymbol{e}_2  & -\Delta & 0 \\
-\mathcal{R}_2\trans\boldsymbol{e}_2  & 0 & -d\Delta
\end{pmatrix}
. 
\]

Problem \eqref{3.4}--\eqref{3.5} is then written as 
\begin{equation}\label{3.6}
\delt\bu+\mathbb{L}_{\mathcal{R}_{1,c}}\bu+\eta \mathbb{P}\b K\b u+\mathbb{P}\b N(\bu)=\b 0. 
\end{equation}
Here $\bu=\trans(\bw,\theta,\psi)$, 
\[
\mathbb{P}=
\begin{pmatrix}
\mathbb{P}_\sigma & \b0 & \b0\\
\trans\b0 & 1 & 0 \\
\trans\b0 & 0 & 1
\end{pmatrix}, 
\ \ \ 
\b K
=
\begin{pmatrix}
\b O & -\Pr\,\b e_2 & \b0 \\
-\trans\b e_2 & 0 & 0 \\
\trans\b0 & 0 & 0 
\end{pmatrix}, 
\]
and 
\[
\b N(\bu) = \b N(\bu,\bu),
\]
where 
\[
\mbox{
$\ds 
\b N(\bu_1,\bu_2)
=
\begin{pmatrix}
\bw_1\cdot\nabla\bw_2 \\
\bw_1\cdot\nabla\theta_2 \\
\bw_1\cdot\nabla\psi_2
\end{pmatrix}
$ 
for $\bu_j
=
\begin{pmatrix}
\bw_j \\
\theta_j \\
\psi_j
\end{pmatrix}
$ 
$(j=1,2)$.}
\] 
We also introduce the adjoint operator $\mathbb{L}_{\mathcal{R}_{1,c}+\eta}^*$ on $\mathbb{X}_\sigma$ which is given by
\[
\mathbb{L}_{\mathcal{R}_{1,c}+\eta}^* = 
\begin{pmatrix}
-\Pr\mathbb{P}_\sigma\Delta & -\Pr(\mathcal{R}_{1,c} + \eta)\mathbb{P}_\sigma\boldsymbol{e}_2 & -\Pr\mathcal{R}_2\mathbb{P}_\sigma\boldsymbol{e}_2\\
-(\mathcal{R}_{1,c} + \eta)\boldsymbol{e}_2  & -\Delta & 0\\
\mathcal{R}_2\boldsymbol{e}_2 & 0 & -d \Delta
\end{pmatrix}
\]
with domain $D(\mathbb{L}_{\mathcal{R}_{1,c}+\eta}^*) = D(\mathbb{L}_{\mathcal{R}_{1,c}+\eta})=\mathbb{Y}_\sigma$. 

The following result on the spectrum of $\mathbb{L}_{\mathcal{R}_{1,c}+\eta}$ was proved by Bona, Hsia, Ma and Wang in \cite{Bona-Hsia-Ma-Wang}.

\begin{prop}\label{prop3.1} {\rm (\cite{Bona-Hsia-Ma-Wang})}  
{\rm (i)} There exist positive numbers $\mathcal{R}_{2*}$ and $\mathcal{R}_2^*$ such that if $\mathcal{R}_2\in[\mathcal{R}_{2*},\mathcal{R}_2^*)$, $\Pr>1$ and $0<d<1$, then the following assertions hold. 
There exist positive constants $\eta_0$, $b_0$ and $\Lambda_0$ such that if $|\eta|\le\eta_0$ then it holds that 
\[
\Sigma\setminus\{\lambda_+(\eta),\lambda_-(\eta)\}\subset\rho(-\mathbb{L}_{\mathcal{R}_{1,c}+\eta}),
\]
where $\Sigma=\{\lambda\in\mathbb{C};\,\re\lambda\ge-b_0|\im\lambda|^2-\Lambda_0\}$; $\lambda_+(\eta)$ and $\lambda_-(\eta)$ are simple eigenvalues of $-\mathbb{L}_{\mathcal{R}_{1,c}+\eta}$ satisfying $\lambda_-(\eta) = \overline{\lambda_+(\eta)}$ and
\[
\lambda_+(0) = ia, \ \ \ \frac{d\re\lambda_+}{d\eta}(0) >0. 
\]
Here $a$ is a positive constant. 

\vspace{1ex} 
{\rm (ii)} Let $\boldsymbol{u}_{\pm}$ be eigenfunctions for the eigenvalues $\pm ia$ of $-\mathbb{L}_{\mathcal{R}_{1,c}}$ and let $\boldsymbol{u}_{\pm}^*$ be eigenfunctions for the eigenvalues $\mp ia$ of the adjoint operator $-\mathbb{L}_{\mathcal{R}_{1,c}}^*$ satisfying $(\boldsymbol{u}_{j}, \boldsymbol{u}_{k}^*) = \delta_{jk}$, where $j, k \in \{+, - \}$. Then $\b{u} = \overline{\b{u}_{+} }$ and $\b{u}^{*} = \overline{\b{u}^{*}_{+} }$; and the eigenprojections $\b P_\pm$ for the eigenvalues $\pm ia$ of $-\mathbb{L}_{\mathcal{R}_{1,c}}$ are given by
\[
\b P_\pm\boldsymbol{u} = (\boldsymbol{u},\boldsymbol{u}_\pm^*)\boldsymbol{u}_\pm.
\]
Furthermore, it holds that  
\[
\frac{d\re\lambda_+}{d\eta}(0) = -\re(\b K\boldsymbol{u}_+,\boldsymbol{u}_+^*)>0.
\]
\end{prop}

We next introduce operators on function spaces of time periodic functions associated with $\mathbb{L}_{\mathcal{R}_{1,c}}$. we define the operators $\mathbb{B}$ and $\mathbb{B}^{*}$ on $L^2(\mathbb{T}_{\frac{2\pi}{a}};\mathbb{X}_\sigma)$ by 
\[
\mathbb{B}=\delt+\mathbb{L}_{\mathcal{R}_{1,c}} \quad  \mbox{ \rm and} \quad  \mathbb{B}^*=-\delt+\mathbb{L}_{\mathcal{R}_{1,c}}^*
\]
with domain $D(\mathbb{B})=D(\mathbb{B}^*)=L^2(\mathbb{T}_{\frac{2\pi}{a}};D(\mathbb{L}_0))\cap H^1(\mathbb{T}_{\frac{2\pi}{a}};\mathbb{X}_\sigma)$. 

We set
\[
\boldsymbol{z}_{\pm} =  e^{\pm iat}\boldsymbol{u}_{\pm},\,\boldsymbol{z}_{\pm}^* =  e^{\pm iat}\boldsymbol{u}_{\pm}^*.
\]
It then follows that $\b z_{\pm}$ and $\b z_{\pm}^*$ are eigenfunctions for the eigenvalue $0$ of $\mathbb{B}$ and $\mathbb{B}^{*}$. Furthermore,  
\[
\langle \boldsymbol{z}_{\pm},\boldsymbol{z}_{\pm}^*\rangle = 1,\,\langle\boldsymbol{z}_{\mp},\boldsymbol{z}_{\pm}^*\rangle = 0.
\]
We define $\hat{\mathcal{P}}_{\pm}$ by
\[
\hat{\mathcal{P}}_\pm\boldsymbol{z} = [\boldsymbol{z}]_\pm\boldsymbol{z}_\pm,\, \mbox{ for }\boldsymbol{z}\in L^2(\mathbb{T}_{\frac{2\pi}{a}};\mathbb{X}_\sigma),
\]
where
\[
[\boldsymbol{z}]_\pm = \langle\boldsymbol{z},\boldsymbol{z}_\pm^*\rangle.
\]
It can be checked that $\hat{\mathcal{P}}_\pm$ are projections satisfying $\hat{\mathcal{P}}_j \hat{\mathcal{P}}_k = \delta_{jk} \hat{ \mathcal{P} }_{j}$ for $j,k\in\{+,-\}$.

\begin{prop}\label{prop3.2}
$0$ is a semisimple eigenvalue of $\mathbb{B}$ and it holds that 
\[
L^2(\mathbb{T}_{\frac{2\pi}{a}};\mathbb{X}_\sigma) = \Ker(\mathbb{B})\oplus R(\mathbb{B}). 
\]
Set  
\[
\hat{\mathcal{P}}_0 = \hat{\mathcal{P}}_+ + \hat{\mathcal{P}}_-, \ \ \ \hat{\mathcal{Q}}_0 = I - \hat{\mathcal{P}}_0. 
\]
Then $\hat{\mathcal{P}}_0$ is an eigenprojection for the eigenvalue $0$ of $\mathbb{B}$ and $\hat{ \mathcal{\mathcal{Q}} }_0$ is a projection on $R(\mathbb{B})$ along $N(\mathbb{B})$. There holds that $\bu \in R(\mathbb{B})$ if and only if $\hat{\mathcal{P}}_0 \bu = 0$, i.e., $[\b z]_+ = [\b z]_- = 0$. 
\end{prop}

We next introduce the associated pressures of the time periodic eigenfunctions. It is known \cite{Galdi, Sohr, Temam3} that there exists the associated pressure $\phi_+\in H_{*,sym}^1(\Omega)$ of $\bu_+=\trans(\bw_+,\theta_+,\psi_+)$, i.e., $u_{+} = \trans(\phi_+, \bu_{+})$ satisfies 
\[
ia 
\begin{pmatrix}
0 \\
\bu_{+} 
\end{pmatrix}
+ L^{\ep}_{\mathcal{R}_{1,c}} u_{+} = 0.  
\]
Here $L^{\ep}_{\mathcal{R}_1}: X \to X$ is the linearized operator around $u = 0$ that is defined by
\begin{align}
& D(L^{\ep}_{\mathcal{R}_1})  =Y, \label{3.8} \\[2ex]
& L^{\ep}_{\mathcal{R}_1}  =
\begin{pmatrix}
0 & \frac{1}{\ep^2}\div & 0 & 0\\
\Pr\nabla & -\Pr\Delta & -\Pr\,\mathcal{R}_{1}\boldsymbol{e}_2 & \Pr\mathcal{R}_2\boldsymbol{e}_2\\
0 & -\mathcal{R}_{1}\trans\boldsymbol{e}_2  & -\Delta & 0\\
0 &-\mathcal{R}_2\trans\boldsymbol{e}_2 & 0 & - d \Delta
\end{pmatrix}.  
\label{3.8'} 
\end{align}
The associated pressure $p_{-}$ of $\bu_{-}=\trans(\bw_-,\theta_-,\psi_-)$ is given by $\phi_{-}=\ov{\phi_+}$, and it holds that $u_{-} = \trans(\phi_{-}, \bu_{-})$ satisfies 
\[
-ia 
\begin{pmatrix}
0 \\
\bu_{-} 
\end{pmatrix}
+ L^{\ep}_{\mathcal{R}_{1,c}} u_{-} = 0.  
\]
We next introduce the adjoint operator $L^{\ep*}_{\mathcal{R}_{1} }: X \to X$ of $L^{\ep}_{\mathcal{R}_{1} }$; it is defined by 
\[
D(L^{\ep*}_{\mathcal{R}_1})=Y,
\]
\[
L^{\ep*}_{\mathcal{R}_1} =
\begin{pmatrix}
0 & - \frac{1}{\ep^2}\div & 0 & 0\\
- \Pr\nabla & -\Pr\Delta & -\Pr\mathcal{R}_{1}\boldsymbol{e}_2 & - \Pr\mathcal{R}_2\boldsymbol{e}_2 \\
0 & -\mathcal{R}_{1}\trans\boldsymbol{e}_2  & -\Delta & 0\\
0 & \mathcal{R}_2\trans\boldsymbol{e}_2 & 0 & -d \Delta
\end{pmatrix}.  
\]
Similarly, we have the associated pressures $\phi_\pm^{*}$ of $\bu_{\pm}^{*}=\trans(\bw_{\pm}^{*},\theta_{\pm}^{*},\psi_{\pm}^{*})$, and it holds that $u_{\pm}^{*} = \trans(\phi_{\pm}^{*}, \bu_{\pm}^{*})$ satisfies 
\[
\mp ia 
\begin{pmatrix}
0 \\
\bu_{\pm}^{*} 
\end{pmatrix}
+ L^{\ep*}_{\mathcal{R}_{1,c}} u_{\pm}^{*} = 0.  
\]

In what follows we set 
\begin{equation}\label{3.9}
u_{\pm} = \trans(\phi_{\pm},\bu_{\pm}), \quad u_{\pm}^{*}=\trans(\phi_{\pm}^{*},\bu_{\pm}^{*}), \quad z_{\pm} = e^{\pm iat} u_{\pm}, \quad z_{\pm}^{*} = e^{\pm iat} u_{\pm}^{*}. 
\end{equation}
We define the operators $\mathcal{P}_{\pm} : \mathcal{X}_{a} \to \mathcal{X}_{a}$ by 
\[
\mathcal{P}_{\pm} u = [\bu ]_{\pm} z_{\pm} 
\]
for $u=\trans(\phi, \bu) \in \mathcal{X}_{a}$, and $\mathcal{P}_{0}, \, \mathcal{Q}_{0} : \mathcal{X}_{a} \to \mathcal{X}_{a}$ by 
\begin{equation}\label{3.10}
\mathcal{P}_{0} = \mathcal{P}_{+} + \mathcal{P}_{-}, \quad \mathcal{Q}_{0} = I - \mathcal{P}_{0}.
\end{equation} 
We note that if $u = \trans(\phi, \bu)$ is real valued, then $\mathcal{P}_{0} u = 2 \re \left( [\bu]_{+} z_{+} \right)$. 

\section{Main results}\label{Main results} 

In this section we state the main results of this paper. 
We fix the parameters $\Pr$, $d$ and $\mathcal{R}_{2}$ in such a way that these parameters satisfies the assumption of Proposition \ref{prop3.1}. 

The artificial compressible system for \eqref{3.4} is written as 
\begin{equation}\label{4.1}
\left\{
\begin{array}{rcl}
\ep^2\delt \phi + \div\bw &=& 0,\\
\delt\bw-\Pr\Delta\bw + \Pr\nabla \phi -\Pr\mathcal{R}_1\theta\boldsymbol{e}_2 + \Pr\mathcal{R}_2\psi\boldsymbol{e}_2 + \bw\cdot\nabla\bw &=& \boldsymbol{0},\\
\delt\theta - \Delta\theta - \mathcal{R}_1\bw\cdot\boldsymbol{e}_2 + \bw\cdot\nabla\theta &=& 0,\\
\delt \psi - d \Delta \psi - \mathcal{R}_2\bw\cdot\boldsymbol{e}_2 + \bw\cdot\nabla \psi &=& 0.
\end{array}
\right.
\end{equation}
Here $\ep$ is a positive parameter, called the artificial Mach number. The system \eqref{4.1} is considered under the boundary condition on the boundary 
$\{x_2=0,1\}$: 
\begin{equation}\label{4.2}
\mbox{$\ds \frac{\del w^1}{\del x_2} = w^2 = \theta = \psi = 0$ on $\{x_2 = 0,1\}$.}
\end{equation}

The linearized problem for \eqref{4.1}--\eqref{4.2} is written as 
\begin{equation}\label{4.3}
\delt u + L^{\ep}_{ \mathcal{R}_{1}} u = 0.
\end{equation}
Here $L^{\ep}_{ \mathcal{R}_{1} }$ is the operator defined in \eqref{3.8}, \eqref{3.8'}. 

We begin with the spectrum of the linearized operator $-L^{\ep}_{\mathcal{R}_1}$ for $\mathcal{R}_{1}$ close to the criticality $\mathcal{R}_{1,c}$ for the incompressible problem. If $\ep \to 0$ in \eqref{4.1}, we formally have the incompressible system. So it could be expected that some part of the spectrum of $-L^{\ep}_{\mathcal{R}_1}$ would be approximated by the one for the incompressible problem although the limiting procedure is a singular limit. In fact, we have the following result.  

\begin{thm}\label{thm4.1} 
{\rm (i)} There exist positive constants $\Lambda_1$, $\ep_1$ and $\eta_1$ such that 
for each $0<\ep\le \ep_1$ there exists a critical value $\mathcal{R}^{\ep}_{1,c}=\mathcal{R}_{1,c}+O(\ep^{2})$ such that  
if $|\eta|\le\eta_1$ with $\eta=\mathcal{R}_{1}-\mathcal{R}^{\ep}_{1,c}$, then
\[
\{\lambda\in\mathbb{C};\,\real\lambda\ge -\Lambda_1\} \setminus \{\lambda^{\ep}_+(\eta),\lambda^{\ep}_-(\eta)\}\subset\rho(-L^{\ep}_{\mathcal{R}_{1,c}^{\ep}+\eta }),
\]
where $\lambda^{\ep}_\pm(\eta)$ are simple eigenvalues of $-L^{\ep}_{\mathcal{R}_{1,c}^{\ep}+\eta}$ satisfying $\lambda^{\ep}_-(\eta)  =  \overline{\lambda^{\ep}_+(\eta)}$,  
\[
\lambda^{\ep}_{\pm}(\eta)  = \lambda_{\pm}(\eta)+O(\ep^2), 
\]
and
\[
\lambda^{\ep}_+(0) = ia^{\ep},
\ \ 
\frac{\partial\re\lambda^{\ep}_+}{\partial\eta}(0) 
= -\re(\b K\boldsymbol{u}_+,\boldsymbol{u}_+^*)+O(\ep^{2})
>0.
\]
Here $a^{\ep}$ is a constant satisfying $a^{\ep}=a+O(\ep^2)$. 

The eigenspaces for $\lambda_\pm^{\ep} (0)$ are spanned by $u^{\ep}_{\pm}$, respectively, where $u^{\ep}_{\pm}$ satisfy $\ov{ u^{\ep}_{-} } = u^{\ep}_{+}$ and $u^{\ep}_{\pm}=u_{\pm}+O(\ep^2)$.

\vspace{1ex}
{\rm (ii)} For $j=\pm$, the eigenprojections $P_j^\ep$ for $\lambda_j^\ep(0)$ $(j=\pm)$ satisfies  
\[
P_j^{\ep}u=(u,u_j^{\ep*})_{\ep} u_j^{\ep},
\]
where $u_{\pm}^{\ep}$ are eigenfunction for the eigenvalues $\lambda_{\pm}^\ep$ of $-L^{\ep}_{\mathcal{R}^{\ep}_{1,c}}$ satisfying $\ov{ u^{\ep}_{-} } = u^{\ep}_{+}$ and $u_{+}^{\ep}=u_{+} + O(\ep^{2})$ in $Y$; $u_{\pm}^{\ep*}$ are the adjoint eigenfunctions for the eigenvalues $\ov{\lambda_j^\ep(0)}$ of $-L^{\ep*}_{\mathcal{R}^{\ep}_{1,c}}$ satisfying $\ov{ u^{\ep*}_{-} } = u^{\ep*}_{+}$, $(u_j^\ep, u_k^{\ep*})_{\ep}=\delta_{jk}$ and $u_{+}^{\ep*}=u_{+}^{*}+O(\ep^{2})$ in $Y$; and $P^{0}_{\pm}$ are projections defined by $P^{0}_{\pm} u = (\bu,\bu_{\pm}^{*}) u_{\pm}$ for $u=\trans(\phi,\bu)$. 
Here $u_{\pm}$ and $u_{\pm}^{*}$ are functions defined in \eqref{3.9}. 
 
Furthermore, if $k\in \mathbb{Z}$ with $k\ge 0$ and $1\le p\le 2$, there exists a positive constant $\tilde{\ep}_{1} = \tilde{\ep}_{1}(k, p, \Omega)$ such that 
\begin{align*}
\|u_{\pm}^{\ep}\|_{H^{k}\times H^{k}} + \|u_{\pm}^{\ep*}\|_{H^{k}\times H^{k}}  & \le C, 
\\ 
\|u_{\pm}^{\ep} - u_{j} \|_{H^{k}\times H^{k}} + \|u_{\pm}^{\ep*} - u_{j}^{*}\|_{H^{k}\times H^{k}}  & \le C \ep^{2}, 
\\ 
\|P_{\pm}^{\ep}u\|_{H^{k}\times H^{k}} & \le C \|u\|_{L^{p}\times L^{p}},
\\
\|(P_{\pm}^{\ep}- P^{0}_{\pm})u\|_{H^{k}\times H^{k}} & \le C\ep^2 \|u\|_{L^{p}\times L^{p}}
\end{align*}
uniformly for $\ep \in (0, \tilde{\ep}_{1}]$.   
\end{thm}

Theorem \ref{thm4.1} can be proved in a similar manner to the arguments in \cite{Kagei-Teramoto}. We here omit the proof. 

In view of Theorem \ref{thm4.1}, it is expected that a Hopf bifurcation will occur when $\mathcal{R}_1$ passes the critical value $\mathcal{R}^{\ep}_{1,c}$. Indeed, the bifurcation of time periodic solutions is proved in the paper \cite{Hsia-Kagei-Nishida-Teramoto2} based on the results in this paper. 

To formulate the time periodic problem for the artificial compressible system \eqref{4.3}, we change the variables $t$ and $u$ into $\tilde{t}$ and $\tilde{u}$, respectively, by $\frac{a^\ep}{a}t = \tilde{t}$ and $u(x,t) = \tilde{u}(x,\tilde{t})$. By omitting the tildes $\tilde{}$ of $\tilde{t}$ and $\tilde{u}$, we deduce from \eqref{4.1} that  
\begin{equation}\label{4.4}
\frac{a^\ep}{a} \delt u + L^{\ep}_{\mathcal{R}_{1,c}^\ep} u + \eta K u + N(u) = 0.
\end{equation} 
Here 
\[
K=
\left(\begin{array}{@{\,}c|ccc@{\,}}
0 & \trans\b0 & 0 & 0 \\ \hline
\b 0 & & &  \\ 
0 & & \b K & \\
0 & & & 
\end{array}\right) 
\] 
and 
\[ 
N(u)= N(u,u),
\]
where 
\[ 
N(u_1,u_2) =
\begin{pmatrix}
0\\
\boldsymbol{N}(\bu_1,\bu_2)
\end{pmatrix}
\]  
for $u_j=\trans(\phi_j,\bu_j)$ $(j=1,2)$. 

The period of a time periodic solution of the nonlinear problem could be expected to be a perturbation of the period of the time periodic eigenfunctions of the linearized problem, so we look for a nontrivial time periodic solution of \eqref{4.4} with period $\frac{2\pi}{a(1+\omega)}$, where $\omega$ is a small number, for sufficiently small $\ep$ and $\eta$. 

We again change the variable $t \mapsto (1+\omega)t$. The problem is then formulated to find a nontrivial time periodic solution with period $\frac{2\pi}{a}$ of the equation 
\begin{equation}\label{4.6}
B^{\ep}u + \frac{a^{\ep}}{a} \omega \delt u+ \eta K u + N (u) = 0, \ \ \ u\in \mathcal{Y}_a.
\end{equation}
Here $B^{\ep}: \mathcal{X}_a \to \mathcal{X}_{a}$ is the operator defined by
\[
B^{\ep} = \frac{a^\ep}{a} \del_t + L^{\ep}_{\mathcal{R}_{1,c}^\ep}
\]
with domain $D(B^{\ep}) = \mathcal{Y}_a$. 

To solve \eqref{4.6} we investigate the spectral properties of $B^{\ep}$. The operator $- \frac{a}{ a^{\ep} (1 + \omega)} L^{\ep}_{ \mathcal{R}_{1,c}^\ep }$ generates an analytic semigroup $e^{-\frac{a t}{ a^{\ep} (1 + \omega)} L^{\ep}_{ \mathcal{R}_{1,c}^\ep } } $ on $X$. See, e.g., \cite{Murata}. Lemma \ref{lem5.1} below shows that $X^{1}$ is invariant under $e^{-\frac{a t}{ a^{\ep} (1 + \omega)} L^{\ep}_{ \mathcal{R}_{1,c}^\ep } } $. We denote the restriction of $e^{-\frac{a t}{ a^{\ep} (1 + \omega)} L^{\ep}_{ \mathcal{R}_{1,c}^\ep } } $ to $X^{1}$ by $\mathscr{V}^{\ep}_{\omega} (t)$: 
\[
\mathscr{V}^{\ep}_{\omega} (t) =  e^{-\frac{a t}{ a^{\ep} (1 + \omega)} L^{\ep}_{ \mathcal{R}_{1,c}^\ep } } \,|_{X^{1} }. 
\]

We set $z_{\pm}^{\ep}=e^{\pm iat}u_\pm^{\ep}$. It is straightforward to verify that $B^{\ep} z_\pm^\ep = 0$. 
We define the operators $\mathcal{P}_{\pm}^{\ep}$ by    
\[
\mathcal{P}_{\pm}^{\ep} u = [u]_{\pm, \ep} z_{\pm}^{\ep}
\]
where
\[
[u]_{\pm, \ep} = \langle u, z_{\pm}^{\ep*} \rangle_\ep 
\]
with $z_{\pm}^{\ep*}=e^{\pm iat}u_\pm^{\ep*}$. 
It then follows that $\mathcal{P}_{\pm}^{\ep}$ are projections on $\mathrm{span}\,\{z_\pm^\ep\}$, respectively. 
As we will see in Theorem \ref{thm5.6} below, $0$ is a semisimple eigenvalue of $B^{\ep}$ and the eigenprojection for the eigenvalue $0$ is given by 
\[
\mathcal{P}^\ep = \mathcal{P}_+^\ep + \mathcal{P}_-^\ep. 
\]

Let $\mathcal{Q}^\ep = I - \mathcal{P}^\ep$. Then $u \in R(\mathcal{Q}^\ep)$ if and only if $[u]_{+,\ep} = [u]_{-,\ep} = 0$. 
Since $[u]_{+,\ep} = \ov{ [u]_{-, \ep} }$ for any real valued function $u$, we see that $u \in R(\mathcal{Q}^\ep)$ if and only if $[u]_{+, \ep} = 0$, when $u$ is a real valued function. We note that $\mathcal{P}_{\pm}^\ep$ can be expressed as 
\[
\mathcal{P}_{\pm}^\ep u = \langle u, \tilde{z}_{\pm}^{\ep *} \rangle_{\ep} \tilde{z}_{\pm}^{\ep}, 
\]
where $\tilde{z}_{\pm}^{\ep} = e^{\pm i a\tau} z_{\pm}^{\ep}$ and $\tilde{z}_{\pm}^{\ep*} = e^{\pm i a\tau} z_{\pm}^{\ep*}$ with a constant $\tau \in [-\frac{\pi}{a}, \frac{\pi}{a})$. 

\begin{thm}\label{thm5.6}
{\rm (i)} If $\ep\in(0,\ep_{1}]$, then
\[
\rho(-B^{\ep})\supset \Sigma_1 \setminus (\cup_{k \in \mathbb{Z}} \{ika^{\ep} \}).
\]
Here $\Sigma_1=\{\lambda\in\mathbb{C};\re\lambda>-\frac{a^{\ep} }{a} \kappa_1 \}$ with $\kappa_1$ being the positive number given in Lemma {\rm \ref{lem5.5}} below; and each $ika^{\ep} $ is a semisimple eigenvalue of $B^{\ep}$ and the corresponding eigenspace is spanned by $ e^{i (1 - k) at}u_{+}^{(\ep)}$ and $ e^{- i (k + 1) at}u_{-}^{(\ep)}$.

\vspace{1ex} 
{\rm (ii)} If $\lambda\in\Sigma_1 \setminus (\cup_{k\in\mathbb{Z}}\{ika^{\ep} \})$, then
\begin{align*}
& (\lambda + B^{\ep})^{-1}F \\
& \quad =  \frac{ 2\pi e^{-\frac{a}{a^{\ep} } \lambda t} }{a^{\ep} ( 1 - e^{-\frac{2\pi}{a^{\ep} }\lambda}) } \left( \left[ e^{-\frac{a}{ a^{\ep} } \lambda (\frac{2\pi}{a}-s) } F \right]_{+, \ep} z_+^\ep + \left[ e^{-\frac{a}{ a^{\ep} } \lambda (\frac{2\pi}{a}-s) } F \right]_{-, \ep} z_-^\ep \right)  
\\
& \quad \quad + \frac{a}{ a^{\ep} } e^{-\frac{a}{ a^{\ep} }\lambda t} \mathscr{V}^{\ep}_{0}(t) \left[ (I- e^{-\frac{2\pi}{ a^{\ep} } \lambda } \mathscr{V}^{\ep}_{0}({\ts\frac{2\pi}{a} }) ) Q^\ep  \right]^{-1} 
\\
& \quad \quad \quad \quad \cdot Q^\ep
 \int_0^{\frac{2\pi}{a}} e^{-\frac{a}{ a^{\ep} } \lambda (\frac{2\pi}{a} - s ) } \mathscr{V}^{\ep}_{0}( {\ts\frac{2\pi}{a} } - s ) F(s) \,ds 
\\
&\quad \quad + \frac{a}{ a^{\ep} } \int_0^t  e^{-\frac{a}{a^{\ep}} \lambda (t-s)} \mathscr{V}^{\ep}_{0}( t - s ) F(s)\,ds.
\end{align*}
Furthermore, 
\[ ||| (\lambda + B^{\ep})^{-1}F |||_{\mathcal{Y}_a} \le \frac{ C}{| 1 - e^{-\frac{2\pi}{ a^{\ep} }\lambda } | } \sum_{j = +, -} \left|\left[ e^{-\frac{a}{a^{\ep} } \lambda (\frac{2\pi}{a}-s) } F \right]_{j,\ep} \right| + C|||F |||_{\mathcal{X}_a}. 
\]

\vspace{1ex}
{\rm (iii)} $0$ is a semisimple eigenvalue of $B^{\ep}$ and 
\[
\mathcal{X}_a = \Ker (B^{\ep}) \oplus R(B^{\ep}), \ \ \ \Ker (B^{\ep}) = \mathrm{span}\,\{z_+^\ep, z_-^\ep\}. 
\]
The projection $\mathcal{P}^\ep$ is an eigenprojection for the eigenvalue $0$ and $\mathcal{Q}^\ep$ is a projection on $R(B^{\ep})$ along $\Ker (B^{\ep})$. 

\vspace{1ex}
{\rm (iv)} Let $F = \trans{(f,\b F)}\in \mathcal{X}_a$. Then the equation 
\[ 
B^{\ep}u = F, \ \  u \in \mathcal{Y}_a
\]
is solvable if and only if $F \in \mathcal{Q}^\ep \mathcal{X}_a$, i.e., $F$ satisfies $[F]_{+,\ep} = [F]_{-,\ep} = 0$. If this condition for $F$ is satisfied, then the problem $B^{\ep}u = F$ with $u \in \mathcal{Q}^\ep \mathcal{Y}_a$ is uniquely solvable and the solution $u \in \mathcal{Q}^\ep \mathcal{Y}_a$ is given by 
\begin{align*}
u(t)  & = \frac{a}{ a^{\ep} } \mathcal{P}^\ep (sF(s))  
\\ 
& \quad + \frac{a}{ a^{\ep} } \mathscr{V}^{\ep}_{0}(t) \left[ (I- \mathscr{V}^{\ep}_{0}({\ts \frac{a}{ a^{\ep} } } ) ) Q^\ep \right] ^{-1}\int_0^{\frac{2\pi}{a}} \mathscr{V}^{\ep}_{0}({\ts \frac{2\pi}{a} }-s ) F(s)\,ds 
\\
& \quad  +\frac{a}{ a^{\ep} } \int_0^t \mathscr{V}^{\ep}_{0}( t - s ) F(s)\,ds. 
\end{align*}
Furthermore, $u$ satisfies
\[
||| u |||_{\ep, \mathcal{Y}_a} \le C \left\{ \int_0^\frac{2\pi}{a} \left( \ep^2\| f \|_2^2 + |\re(\b F, \bu)| +\ep^2 \|\b F\|_2^2 + \ep^6 \|\delx f\|_2^2 \right) \, dt \right\}^{\frac{1}{2}}
\]
and 
\[
||| u |||_{\ep, \mathcal{Y}_a} \le C ||| F |||_{\ep, \mathcal{X}_a} 
\]
uniformly in $\ep\in(0,\ep_{1}]$.
\end{thm}

We see from Theorem \ref{thm5.6} that $0$ is a semisimple eigenvalue of $B^{\ep}$ whose eigenspace is spanned by $z_+^\ep$ and $z_-^\ep$. 
To prove the Hopf bifurcation in \eqref{4.6}, we will employ the Lyapunov-Schmidt method, i.e., we decompose \eqref{4.6} to the $\mathcal{P}^\ep$ part (finite dimensional part) and $\mathcal{Q}^\ep$ part (infinite dimensional part). 
In a standard setting, the term $\omega \delt u$ on the left-hand side of the $\mathcal{Q}^\ep$ part of \eqref{4.6} is regarded as a perturbation of the term $B^{\ep} u$ since $B^{\ep}$ has a bounded inverse on the $\mathcal{Q}^\ep$ part and $\omega$ is sufficiently small. However, in the setting of this paper, if we do so, we cannot establish uniform estimates with respect to $\ep$ due to the time derivative of the pressure $\delt \phi$. (See the estimate in Lemma \ref{lem5.1} below.) To avoid this kind of ``$\ep$ loss'', we will put $\omega \delt u$ in the principal part. 
For this purpose we prepare the following lemma which is analogous to Theorem \ref{thm5.6}. 

We define the operator $B^{\ep}(\omega)$ on $\mathcal{X}_a$ by 
\[
D(B^{\ep}( \omega))=\mathcal{Y}_a, \ \ B^{\ep}( \omega) u = \frac{a^{\ep}}{a} (1+\omega)\delt u + L^{\ep}_{\mathcal{R}_{1,c}^\ep} u \ \ (u \in D(B^{\ep}( \omega))). 
\]
  
\begin{lem}\label{lem5.8} 
There exist a positive constant $\ep_{1}$ such that the following assertion holds.  
For a given $F = \trans(f, \b F) \in \mathcal{Q}^\ep \mathcal{X}_a$, there exists a unique solution $u \in \mathcal{Q}^\ep \mathcal{Y}_a$ of 
$B^{\ep}( \omega) u = F $, and the solution $u$ satisfies the estimate 
\[
||| u |||_{\ep, \mathcal{Y}_a} \le C \left\{ \int_0^\frac{2\pi}{a} \left( \ep^2\| f \|_2^2 + |\re(\b F, \bu)| +\ep^2 \|\b F\|_2^2 + \ep^6 \|\delx f\|_2^2 \right) \, dt \right\}^{\frac{1}{2}}
\]
and  
\[
||| u |||_{\ep, \mathcal{Y}_a} \le C||| F |||_{\ep, \mathcal{X}_a}
\]
uniformly for $0 <\ep \le \ep_{1}$ and $|\omega| \le \frac{1}{4}$.  
\end{lem}

\section{Proofs of Theorem \ref{thm5.6} and Lemma \ref{lem5.8} }\label{Spectrum of B}
 
In this section we investigate the spectrum of $B^{\ep}$ near the imaginary axis. We shall prove Theorem \ref{thm5.6} and Lemma \ref{lem5.8}. In what follows we assume that $\mathcal{R}_{1}$ is in the interval $[\frac{1}{2}\mathcal{R}_{1,c}, \frac{3}{2}\mathcal{R}_{1,c}]$. 

We first show that $ika\,(k\in\mathbb{Z})$ are semisimple eigenvalues of $B^{\ep}$. In fact, it is easy to verify that $B^{\ep}( e^{i(k+1)at}u_{+}^{(\ep)}) = ika e^{i(k+1)at}u_{+}^{(\ep)}$. Therefore, $ika\,(k\in\mathbb{Z})$ are eigenvalues. To show that $ika\,(k\in\mathbb{Z})$ are semisimple eigenvalues, we prepare the following lemma.

\begin{lem}\label{lem5.1}
Let $\beta$ and $T$ be positive constants. Suppose that $u_0\in X^1$ and let $F=\trans(f,\b g,h,k)\in  L^2(0,T; X)$. 
Then there is a unique solution $u(t)\in L^2(0, T;Y) \cap H^1(0, T;X)$ of the problem
\begin{align}\label{5.1}
\beta \del_t u + L^{\ep}_{\mathcal{R}_{1,c}^\ep}u & =  F, \\ \label{5.2}
u(0)  & = u_0
\end{align}
with estimate
\begin{align*}
& \beta ||| u(t) |||_{\ep, X^1}^2 
+ \int_0^t \left(|||\delx u |||_\ep^2 + \beta^{2} \ep ^2 |||\delt u |||_\ep^2 + \ep^2 \|\delx ^2\bu \|_2^2 + \beta^{2} \ep^6 \|\delx\delt \phi\|_2^2 \right) \,ds \\
& \ \le C e^{\frac{C}{\beta} t} \left\{ ||| u_0 |||_{\ep, X^1}^2  
+ \int_0^t \left( \ep^2\| f \|_2^2 + |\re(\b F, \bu)| +\ep^2 \|\b F\|_2^2 + \ep^6 \|\delx f\|_2^2 \right) \, ds \right\}
\end{align*}
uniformly for $0 \le t \le T$ and $0 <\ep \le \ep_{1}$ and $\beta$. 
Furthermore, if $\beta \ge \beta_{1}$, then the estimate 
\[
||| u |||_{\ep, \mathcal{Y}(0, T) } \le C \left\{ ||| u_0 |||_{\ep, X^1} + ||| F |||_{\ep, \mathcal{X}(0, T) } \right \}
\]
holds uniformly for $0 <\ep \le \ep_{1}$ with a constant depending only on $\beta_{1}$. 
\end{lem}

The existence of a solution of \eqref{5.1}--\eqref{5.2} is proved by a standard iteration argument based on the parabolic theory.  The estimates in Lemma 5.1 follow from the following basic estimates. We set 
\[
D(\bu)=\|\nabla\bw\|_2^2 + \|\nabla\theta\|_2^2 + \tau\|\nabla\psi\|_2^2 
\] 
for $\bu=\trans(\bw,\theta,\psi)$. 

\begin{lem}\label{lem5.2}
There are positive constants $c$ and $C$ such that the following estimates hold for a solution $u=\trans(\phi,\bu)\in L^2(0,\frac{2\pi}{a}; Y) \cap H^1(0,\frac{2\pi}{a}; X)$ of \eqref{5.1}--\eqref{5.2} uniformly in $0< \ep \le 1$ and $\beta$:  
\begin{equation} \label{5.3}
\beta \frac{d}{dt} |||u|||_\ep^2 +D(\bu)  \le 2\re(F,u)_\ep+C\left\{ \|\theta\|_2^2 + \|\psi\|_2^2\right\}, 
\end{equation}
\begin{equation}\label{5.4}
\begin{split}
& \beta \frac{d}{dt} |||\del_{x_1} u|||_\ep^2 +D(\del_{x_1} \bu)+c\beta^{2}\ep^4\|\delt\del_{x_1}\phi\|_2^2 
\\  
& \quad \le \ep^2\re(\del_{x_1}f,\del_{x_1}\phi)_{L^2} + C\left\{ \ep^4\|\del_{x_1}f\|_2^2 + \|\b F \|_2^2 + \|\bu\|_2^2 \right\},  
\end{split}
\end{equation}
\begin{equation}\label{5.5}
\begin{split}
& \beta \frac{d}{dt} \left(D(\bu)-2\re(\phi,\div\bw)_{L^2}\right) + \beta^{2} |||\delt u|||_\ep^2 
\\  
& \quad \le C\left\{ \frac{1}{\ep^2}\|\div\bw\|_2^2 + |||F|||_\ep^2 + \|\bu\|_2^2 \right\},  
\end{split}
\end{equation}
\begin{equation}\label{5.6}
\begin{split}
& \beta \frac{d}{dt} \ep^2\|\del_{x_2}\phi\|_2^2 + \|\del_{x_2}\phi\|_2^2 + c \beta^{2} \ep^4\|\delt \del_{x_2}\phi\|_2^2 
\\  
& \quad \le C\left\{\ep^4\|\del_{x_2} f \|_2^2 + \|\b g\|_2^2 + D(\del_{x_1}\bu) + \beta^{2} \|\delt \bw\|_2^2 + \|\bu\|_2^2 \right\},  
\end{split}
\end{equation}
and 
\begin{equation}\label{5.7}
\begin{split}
& \|\delx \phi\|_2^2 + \|\delx^2 \bu\|_2^2  \\ 
& \quad \le  C\left\{\ep^4\| f \|_{H^1}^2 + \|\b F\|_2^2 + \beta^{2} \ep^4\|\delt \phi\|_{H^1}^2  + \beta^{2} \|\delt \bu\|_2^2 + \|\bu\|_2^2 \right\}.  
\end{split}
\end{equation}
\end{lem}

\vspace{1ex}
\noindent
{\bf Proof.} We write \eqref{5.1}--\eqref{5.2} in the form: 
\begin{align}
\label{5.8}
\beta \delt u + A^{\ep} u & =\tilde F,  \\
\label{5.9}
u(0) & = u_0. 
\end{align}
Here $A^{\ep}$ is the operator on $X$ defined by  
\[
D(A^{\ep})=Y, \ \ \ 
A^{\ep} =
\begin{pmatrix}
0 & \frac{1}{\ep^2}\div & 0 & 0 \\
\Pr\nabla & -\Pr\Delta & 0 & 0 \\
0 & \trans\b0 & -\Delta & 0 \\
0 & \trans\b0 &  0 & -\tau \Delta 
\end{pmatrix}, 
\]
and $\tilde{F}$ is given by 
\[
\tilde F=
\begin{pmatrix}
f \\ 
\tilde{\b g} \\
\tilde h \\
\tilde k
\end{pmatrix}
=F+(A^{\ep} - L^{\ep}_{\mathcal{R}_{1,c}^\ep } ) u
=
\begin{pmatrix}
f \\ 
\b{g} \\
h \\
k
\end{pmatrix}
+
\begin{pmatrix}
0 \\ 
\Pr (\mathcal{R}^{\ep}_{1, c}\theta -\mathcal{R}_2\psi)\b e_2 \\
\mathcal{R}^{\ep}_{1, c} \bw\cdot \b e_2\\
\mathcal{R}_2 \bw\cdot \b e_2
\end{pmatrix}. 
\]

Taking the inner product of \eqref{5.8} with $u$, we have 
\begin{equation}\label{5.10}
\frac{\beta}{2}\frac{d}{dt} |||u|||_\ep^2 + D(\bu) + (\div\bw, \phi)_{L^2} + (\nabla \phi, \bw)_{L^2} 
= (\tilde F, u)_\ep. 
\end{equation} 
By integration by parts, we have $(\div\bw,\phi)_{L^2} + (\nabla \phi, \bw)_{L^2} = 2i\im(\div\bw, \phi)_{L^2}$. By the Poincar\'e inequality, we obtain $\re(\tilde F, u)_\ep \le \re(F, u)_\ep + \frac{1}{2}\|\nabla \bw\|_2^2+C\left\{ \|\theta\|_2^2+\|\psi\|_2^2\right\}$. Therefore, the real part of \eqref{5.10} gives the estimate \eqref{5.3}. 

Applying $\del_{x_1}$ to \eqref{5.8}, we similarly have 
\begin{align*}
& \frac{\beta}{2}\frac{d}{dt} |||\del_{x_1} u|||_\ep^2 + D(\del_{x_1} \bu)  =  \re(\del_{x_1}\tilde F,\del_{x_1}u)_\ep \\
& \quad \le \frac{1}{4}D( \del_{ x_{1} } \bu) + \ep^2\re(\del_{x_1}f,\del_{x_1}\phi)_{L^2} + C\left\{ \|\b F\|_2^2+\|\bu\|_2^2 \right\}, 
\end{align*}
and hence, 
\begin{equation}\label{5.11}
\begin{split}
& \beta \frac{d}{dt} |||\del_{x_1} u|||_\ep^2 + \frac{3}{2}D(\del_{x_1} \bu) \\
& \quad  \le \ep^2\re(\del_{x_1}f,\del_{x_1}\phi)_{L^2} + C\left\{ \|\b F\|_2^2+\|\bu\|_2^2 \right\}. 
\end{split}
\end{equation}
On the other hand, the first row of the equation \eqref{5.8} is written as 
\begin{equation}\label{5.12}
\beta \ep^2\delt \phi + \div\bw = \ep^2 f,
\end{equation}
which yields   
\[
\beta^{2} \ep^4\|\delt \del_{x_1}\phi \|_2^2 \le C\left\{ D(\del_{x_1}\bu)+\ep^4\|\del_{x_1} f \|_2^2 \right\}.
\]
Combining this with \eqref{5.11}, we obtain the estimate \eqref{5.4}. 
 
As for the estimate \eqref{5.5}, we take the inner product of \eqref{5.8} with $\delt u$ to obtain 
\[
\beta |||\delt u|||_\ep^2 + \frac{1}{2}\frac{d}{dt} D(\bu) + (\nabla \phi, \delt\bw)_{L^2} +(\div \bw, \delt\phi)_{L^2} 
=(\tilde F, \delt u)_\ep. 
\]
Taking the real part of this equation, we have 
\begin{equation}\label{5.13} 
\begin{split}
& \frac{\beta}{2}\frac{d}{dt} D(\bu) + |||\delt u|||_\ep^2 + \re\left\{ (\nabla \phi, \delt\bw)_{L^2} +(\div \bw, \delt\phi)_{L^2} \right\} 
\\
& \quad \le \frac{\beta}{2}|||\delt u|||_\ep^2 + \frac{C}{\beta} \left\{|||F|||_\ep^2+\|\bu\|_2^2\right\}. 
\end{split}
\end{equation}
We consider the third term on the left-hand side of \eqref{5.13}. 
By \eqref{5.12}, we have 
\begin{align*}
& (\div \bw,\delt \phi)_{L^2} + (\nabla \phi, \delt\bw)_{L^2} \\
& \quad = -\frac{d}{dt}(\phi, \div\bw)_{L^2} + 2\re(\div\bw,\delt \phi)_{L^2} \\
& \quad = -\frac{d}{dt}(\phi,\div\bw)_{L^2} -\frac{2}{\beta \ep^2}\|\div \bw\|_2^2 + \frac{2}{\beta} \re(\div\bw, f)_{L^2}. 
\end{align*}
This, together with \eqref{5.13}, implies  
\begin{align*}
& \frac{d}{dt} \left( D(\bu) - 2\re(\phi,\div\bw)_{L^2}\right) + \beta |||\delt u|||_\ep^2 
\\
& \quad \le  \frac{4}{\beta \ep^2}\|\div \bw\|_2^2 - \frac{4}{\beta} \re(\div\bw, f)_{L^2}+ \frac{C}{\beta} \left\{|||F|||_\ep^2+\|\bu\|_2^2\right\}
\\
& \quad \le \frac{C}{\beta} \left\{\frac{1}{\ep^2}\|\div \bw\|_2^2 + |||F|||_\ep^2+\|\bu\|_2^2 \right\}.   
\end{align*}
We thus obtain the estimate \eqref{5.5}. 
 
To obtain the estimate \eqref{5.6}, we make use of \eqref{5.12} and the third row of the equation \eqref{5.8} which takes the form: 
\begin{equation}\label{5.14} 
\beta \Pr^{-1}\delt w^2 - \Delta w^2 +\del_{x_2}\phi = \tilde g^2, 
\end{equation}
where $\tilde g^2 = \Pr^{-1} g^2 +(\mathcal{R}^{\ep}_{1,c}\theta - \mathcal{R}_2\psi)$. 

We compute $\del_{x_2}\eqref{5.12} + \eqref{5.14}$ to obtain 
\begin{equation}\label{5.15}
\beta \ep^2\delt\del_{x_2} \phi  + \del_{x_2} \phi  = \ep^2\del_{x_2} f + \tilde{g}^2 -\{\beta \Pr^{-1}\delt w^2 + \del_{x_1}\del_{x_2}w^1 - \del_{x_1}^2 w^2 \}.
\end{equation}
We take the $L^2$ inner product of \eqref{5.15} with $\del_{x_2}\phi$ and then take the real part of the resulting equation to obtain 
\begin{equation*}
\begin{split}
& \frac{\beta}{2}\frac{d}{dt}\left(\ep^2\|\del_{x_2} \phi \|_2^2\right) + \|\del_{x_2} \phi \|_2^2 \\ 
& \quad  = \re \left(\ep^2\del_{x_2} f + \tilde{g}^2-\left\{ \beta \Pr^{-1}\delt w^2 + \del_{x_1}\del_{x_2}w^1 + \del_{x_1}^2w^2 \right\}, \del_{x_2} \phi \right)_{L^2} \\
& \quad \le \frac{1}{4}\|\del_{x_2} \phi \|_2^2 + C\left\{\ep^4\|\del_{x_2} f \|_2^2 + \|\b g\|_2^2 + D(\del_{x_1}\bu) + \beta^{2} \|\delt \bw\|_2^2 + \|\bu\|_2^2 \right\}, 
\end{split}
\end{equation*}
and hence, 
\begin{equation}\label{5.16}
\begin{split}
& \beta \frac{d}{dt}\left(\ep^2\|\del_{x_2} \phi \|_2^2\right) + \frac{3}{2}\|\del_{x_2} \phi \|_2^2 \\ 
& \quad \le C\left\{\ep^4\|\del_{x_2} f \|_2^2 + \|\b g\|_2^2 + D(\del_{x_1}\bu) + \beta^{2} \|\delt \bw\|_2^2 + \|\bu\|_2^2 \right\}.
\end{split}
\end{equation}
We also see from \eqref{5.15} that 
\begin{align*}
\beta^{2} \ep^4\|\delt\del_{x_2}\phi\|_2^2  \le & C\left\{ \|\del_{x_2} \phi \|_2^2+\ep^4\|\del_{x_2} f \|_2^2 + \|\b g\|_2^2 \right. \\ 
 & \quad \, \left. + D(\del_{x_1}\bu) + \beta^{2} \|\delt \bw\|_2^2 + \|\bu\|_2^2 \right\}.
\end{align*}
This, together with \eqref{5.16}, give the estimate \eqref{5.6}. 

To obtain the estimate \eqref{5.7}, we make use of the elliptic estimate for $-\Delta$ under the Dirichlet boundary condition and the following estimate for the Stokes system. 

\begin{lem}\label{stokes}
If $f \in H^1_{*,sym}(\Omega)$ and $\b g \in \b L_{sym}^2(\Omega)$, then there exists a unique pair of functions $\trans(p,\bv) \in H^1_{*,sym}(\Omega)\times [\b H_b^2(\Omega)\cap \b L_{sym}^2(\Omega)]$ satisfying  
\[
\left\{
\begin{array}{rcl}
\div \bv &=& f\\
-\Delta\bv + \nabla p &=& \b g 
\end{array}
\right. \mbox{ \rm in } \Omega.
\]
Furthermore, $\trans(p,\bv)$ satisfies the following estimates:
\[
\|p\|_2 + \|\nabla\bv\|_2 \le C\{\|f\|_2 + \|\b g\|_{H^{-1}}\},
\]
\[
\|\nabla p\|_2 + \|\nabla^2\bv\|_2 \le C\{\|f\|_{H^1} + \|\b g\|_2\}.
\]
\end{lem}

Lemma \ref{stokes} can be proved by using Fourier series expansions. We omit the proof. (See \cite{Galdi, Sohr, Temam3} for the case of the non-slip boundary condition.)  

Since 
\[
\left\{
\begin{array}{rcl}
\div \bw & = & \ep^2 f - \beta \ep^2\delt\phi, \\
-\Delta \b w + \nabla \phi & = & \Pr^{-1}\left( \tilde{\b g} - \beta \delt \bw \right),  \\
-\Delta \theta & = & \tilde h - \beta \delt \theta, \\ 
-\tau \Delta \psi & = & \tilde k - \beta \delt \psi, 
\end{array}
\right. 
\]
and $\bu = \trans( \bw, \theta, \psi) \in \b Y$, we obtain the estimate \eqref{5.7}. 
This completes the proof of Lemma \ref{lem5.2}. 
\hfill$\square$.  

\vspace{1ex}
\noindent
{\bf Proof of Lemma \ref{lem5.1}.} 
We first make an observation on the quantity $D(\bu)-2\re(\phi,\div\bw)_{L^2}$ which appears in the estimate \eqref{5.5}. 
We set 
\[
E_1(u)=\beta |||u|||_\ep^2+c_1\ep^2\left\{D(\bu)-2\re(\phi,\div\bw)_{L^2}\right\},
\]
where $c_1$ is a positive constant (independent of $\ep$ and $\beta$) to be determined later. 
Since 
\[
\frac{1}{2}D(\bu) - 8\|\phi\|_2^2 \le D(\bu)-2\re(\phi,\div\bw)_{L^2} \le \frac{3}{2}D(\bu)  + 8\|\phi\|_2^2, 
\]
we see that if $c_1= \frac{\beta}{16}$, then it holds that 
\begin{equation}\label{5.17}
C^{-1}\beta \left\{|||u|||_\ep^2+\ep^2 D(\bu)\right\} \le E_1(u) \le C \beta \left\{|||u|||_\ep^2+\ep^2 D(\bu)\right\}
\end{equation}
uniformly for $\ep > 0$. 

We compute $\eqref{5.3} + \ep^2 c_{0} \times \eqref{5.4} + c_1\ep^2\times \eqref{5.5} + c_2\ep^2\times \eqref{5.6} + c_3\ep^2\times \eqref{5.7}$. 
In view of \eqref{5.17} and the Poincar\'e inequality, one can take $c_j$ $(j=0,1,2,3)$ suitably small (but independent of $\ep$) to arrive at the estimate
\begin{equation}\label{5.18}
\begin{split}
& \frac{d}{dt}E(u) + |||\delx u|||_\ep^2 + \beta^{2} \ep ^2 |||\delt u|||_\ep^2 + \ep^2 \|\delx ^2\bu \|_2^2 + \beta^{2} \ep^6 \|\delx\delt \phi\|_2^2  \\
& \quad  \le C \left\{ |\re(F,u)_\ep| + \ep^4 |\re(\del_{x_1}f,\del_{x_1}\phi)_{L^2}|  \right. \\ 
& \quad \quad \quad \left. + \ep^2 ||| F |||_\ep^2 + \ep^6 \|\delx f\|_2^2 + \|\theta \|_2^2 + \|\psi \|_2^2 \right\} 
\end{split}
\end{equation} 
uniformly for $0 \le t \le \frac{2\pi}{a}$ and $0 <\ep \le 1$. 
Here $E(u)$ is an energy functional equivalent to $\beta \left\{ ||| u |||_\ep^2 + \ep^2 |||\delx u|||_\ep^2 \right\}$, i.e.. 
\[
C^{-1} \beta \left\{||| u |||_\ep^2 + \ep^2 |||\delx u|||_\ep^2\right\} \le E(u) \le C \beta \left\{||| u |||_\ep^2 + \ep^2 |||\delx u|||_\ep^2\right\}
\]
uniformly for $0 <\ep \le 1$ and $\beta$. 

In fact, by the Poincar\'e inequality $\|\bu\|_2^2 \le CD(\bu)$ and \eqref{5.17}, we can find positive constants $ c_{0}$ and $c_{1}$ such that $C\left\{ (c_{0} + c_{1} \ep^2)\|\bu\|_2^2+  c_{1} \|\div\bw\|_2^2 \right\}$ on the right-hand side of $\eqref{5.3} + c_{0} \ep^2\times \eqref{5.4} + c_1\ep^2\times \eqref{5.5}$ is absorbed into the left-hand side, and obtain 
\begin{equation}\label{5.19}
\begin{split}
& \frac{d}{dt}\left(E_1(u) + \beta |||\del_{x_1} u|||_\ep^2 \right) + D_1(u) 
\\
& \quad \le C \left\{ |\re(F,u)_\ep| + \ep^4 |\re(\del_{x_1}f,\del_{x_1}\phi)_{L^2}|  \right. \\
&  \quad \quad \quad \left. + \ep^2 ||| F |||_\ep^2 + \ep^6 \|\del_{x_1}  f\|_2^2 + \|\theta \|_2^2 + \|\psi \|_2^2 \right\},   
\end{split}
\end{equation}
where $D_1(u)=D(\bu) + \ep^2D(\del_{x_1}\bu) + \beta^{2} \ep^2||| \delt u |||_\ep^2 + \beta^{2} \ep^6\|\delt\del_{x_1}\phi\|_2^2$. 

We can then find a positive constant $c_2$ so that $c_2\ep^2C\{ D(\del_{x_1}\bu) + \beta^{2} \|\delt \bw\|_2^2 + \|\bu\|_2^2 \}$ on the right-hand side of $\eqref{5.19} + c_2\ep^2\times \eqref{5.6}$ is absorbed into the left-hand side to obtain 
\begin{equation}\label{5.20}
\begin{split}
\frac{d}{dt}E(u) + D_2(u) & \le C \left\{ |\re(F,u)_\ep| + \ep^4 |\re(\del_{x_1}f,\del_{x_1}\phi)_{L^2}|  \right. \\
&  \quad \quad \left. + \ep^2 ||| F |||_\ep^2 + \ep^6 \|\delx  f\|_2^2 + \|\theta \|_2^2 + \|\psi \|_2^2 \right\},   
\end{split}
\end{equation}
where $E(u)=E_1(u) + \beta |||\del_{x_1} u|||_\ep^2 + c_2\beta\ep^4\|\del_{x_2} \phi\|_2^2$ and $D_2(u)= D(\bu) + \ep^2D(\del_{x_1}\bu) + \beta^{2} \ep^2||| \delt u |||_\ep^2 + \ep^2 \|\del_{x_2}\phi\|_2^2 + \beta^{2} \ep^6\|\delt\delx \phi\|_2^2$. 
By \eqref{5.17}, we see that $E(u)$ is an energy functional equivalent to $\beta \{ ||| u |||_\ep^2 + \ep^2 |||\delx u|||_\ep^2 \}$.     

Similarly, we can find a positive constant $c_3$ so that $c_3\ep^2C \{ \beta^{2} \ep^4\|\delt \phi\|_{H^1}^2 + \beta^{2} \|\delt \bu\|_2^2 + \|\bu\|_2^2 \}$ on the right-hand side of $\eqref{5.20} + c_3\ep^2\times \eqref{5.7}$ is absorbed into the left-hand side, and we arrive at the estimate \eqref{5.18}. 

Since 
\begin{align*}
|\re(F,u)_\ep| & \le \kappa\ep^2\|\phi\|_2^2+\frac{\ep^2}{\kappa}\|f\|_2^2 + \kappa \|\delx \bu\|_2^2+\frac{1}{\kappa}\|\b F \|_{(\b X^1)*}^2 \\
& \le \kappa\ep^2\|\delx \phi\|_2^2+\frac{\ep^2}{\kappa}\|f\|_2^2 + \kappa \|\delx \bu\|_2^2+\frac{1}{\kappa}\|\b F \|_{(\b X^1)*}^2 
\end{align*}
and 
\[
\ep^4|\re(\del_{x_1}f,\del_{x_1}\phi)_{L^2}| \le \kappa\ep^2\|\del_{x_1} \phi\|_2^2+\frac{\ep^6}{\kappa}\|\del_{x_1} f\|_2^2
\]
for any $\kappa>0$, by taking $\kappa = \frac{1}{6C}$, we obtain 
\begin{equation}\label{5.21}
\begin{split}
& \frac{d}{dt}E(u) + |||\delx u|||_\ep^2 + \beta^{2} \ep ^2 |||\delt u|||_\ep^2 + \ep^2 \|\delx ^2\bu \|_2^2 + \beta^{2} \ep^6 \|\delx\delt \phi\|_2^2  \\
& \quad  \le C \left\{ \ep^2\| f \|_2^2 + \|\b F \|_{(\b X^1)*}^2 + \ep^2 \|\b F\|_2^2 + \ep^6 \|\delx f\|_2^2 + \|\theta \|_2^2 + \|\psi \|_2^2 \right\} 
\end{split}
\end{equation} 
uniformly for $0 \le t \le T$ and $0 <\ep \le \ep_{1}$. 
Since $\|\theta \|_2^2 + \|\psi \|_2^2\le \|\bu\|_2^2$, by using the Gronwall inequality, we obtain the desired estimate. 
This completes the proof. 
\hfill$\square$

\vspace{1ex}
We next investigate the spectrum of the operator $e^{-\frac{a t}{ a^{\ep} (1 + \omega)} L^{\ep}_{ \mathcal{R}_{1,c}^\ep } }$ with $t = \frac{2\pi}{a}$ on the space $X^1$. Recall that the restriction of $e^{-\frac{a t}{ a^{\ep} (1 + \omega)} L^{\ep}_{ \mathcal{R}_{1,c}^\ep } } $ to $X^{1}$ is denoted by $\mathscr{V}^{\ep}_{\omega} (t)$: 
\[
\mathscr{V}^{\ep}_{\omega} (t) =  e^{-\frac{a t}{ a^{\ep} (1 + \omega)} L^{\ep}_{ \mathcal{R}_{1,c}^\ep } } \,|_{X^{1} }. 
\]

We define the projections $P^{\ep}$ and $Q^\ep$ by  
\[
\mbox{ \rm $P^{\ep} = P^{\ep}_{+} + P^{\ep}_{-}$ and $Q^\ep = I - P^\ep$, } .
\] 
respectively. 
Since $\pm ia^{\ep}$ are simple eigenvalues of $ - L^{\ep}_{\mathcal{R}^{\ep}_{1,c}}$, we have   
\[
P^{\ep}_{\pm} \mathscr{V}^{\ep}_{\omega}(t) = \mathscr{V}^{\ep}_{\omega}(t) P^{\ep}_{\pm}  = e^{ \pm \frac{ iat }{ 1 + \omega} } P^{\ep}_{\pm}. 
\]

We shall derive the decay estimate for the $Q^{\ep}$ part of $\mathscr{V}^{\ep}_{\omega}(t)$  in Lemma \ref{lem5.5} below.   

In the remaining of this section we will use the same letters $\ep_{1}$, $\delta_{1}$ and $\omega_{1}$ for the ranges of $\ep$, $\delta$ and $\omega$ if $\ep_{1}$, $\delta_{1}$ and $\omega_{1}$ depend only on $\Pr$, $\mathcal{R}_{2}$, $\tau$, $\mathcal{R}_{1,c}$ and $\Omega$.  

\begin{lem}\label{lem5.5}
There exist positive constants $\ep_{1}$ and $\kappa_1$ such that the estimate 
\[
|||\mathscr{V}^{\ep}_{0}(t) Q^\ep u_0 |||_{\ep, X^1} \le Ce^{-\kappa_1 t} |||u_0|||_{\ep, X^1}
\]
holds uniformly in $t\ge0$ and $0<\ep\le\ep_{1}$. 
Furthermore, if $\re\lambda > -\kappa_1$, then $(I - e^{-\frac{2\pi}{a}\lambda} \mathscr{V}^{\ep}_{0} (\frac{2\pi}{a}) ) Q^\ep$ has a bounded inverse on $Q^\ep X^1$ and its inverse $[(I - e^{-\frac{2\pi}{a}\lambda} \mathscr{V}^{\ep}_{0} (\frac{2\pi}{a}) ) Q^\ep]^{-1}$ satisfies
\[
\| [ (I - e^{-\frac{2\pi}{a}\lambda} \mathscr{V}^{\ep}_{0} (\ts \frac{2\pi}{a} ) ) Q^\ep ]^{-1} \|
\le \frac{C}{1-e^{-\frac{2\pi}{a}(\re\lambda + \kappa_1)}}
\]
uniformly in $\ep\in(0,\ep_{1}]$.
\end{lem}

\vspace{1ex}
\noindent 
{\bf Proof.} 
We consider the problem \eqref{5.1}--\eqref{5.2} with $F=0$ and $\beta = \frac{ a^{\ep} }{ a}$. 
By the Poincar\'e inequality, we see that there exists a positive constant $\kappa_0$ such that 
$|||\delx u|||_\ep^2 \ge 2\kappa E(u)$ for $0\le \kappa \le \kappa_0$ and $0 < \ep \le 1$. 
It then follows from \eqref{5.21} that if $0 < \ep \le \ep_{1}$, then the solution $u=\trans(\phi,\bw,\theta,\psi)$ of \eqref{5.1}--\eqref{5.2} with $F=0$ and $\frac{ a^{\ep} }{a}$ satisfies  
\[
\frac{d}{dt}E(u) + 2\kappa E(u) \le C\left\{ \|\theta\|_2^2 + \|\psi\|_2^2 \right\}.
\]
We thus obtain 
\begin{equation}\label{5.22} 
E(u(t)) \le e^{-2\kappa t} E(u_0) + C e^{-2\kappa t} \int_0^t e^{2\kappa s} \left\{ \|\theta\|_2^2 + \|\psi\|_2^2 \right\}\, ds.
\end{equation} 

Let us prove that there exists a positive number $\kappa$ such that 
\begin{equation}\label{5.23}
\int_0^t e^{2\kappa s} \left\{ \|\theta\|_2^2 + \|\psi\|_2^2 \right\}\, ds \le C E(u_0)
\end{equation} 
uniformly in $t \ge 0$ and $0< \ep \le \ep_{1}$ for some positive constant $\ep_{1}$. 

To prove \eqref{5.23}, we first observe that 
\[
u(t) =\frac{1}{2\pi i}\int_{\Gamma} e^{\lambda t}(\lambda + {\ts \frac{a}{ a^{\ep} } } L^{\ep}_{\mathcal{R}_{1,c}^\ep} )^{-1} u_0 \, d\lambda,
\]
where $\Gamma = \{\lambda \in \mathbb{C} ; \re\lambda\ge -\ep \tilde{b}_0 |\im\lambda| + \tilde\Lambda_0\}$ for some positive constants $\tilde{b}_0$ and $\tilde\Lambda_0$. 
We set 
\[
\trans(\phi_\lambda[u_0], \bw_\lambda[u_0], \theta_\lambda[u_0], \psi_\lambda[u_0])=(\lambda + {\ts \frac{a}{ a^{\ep} } } L^{\ep}_{\mathcal{R}_{1,c}^\ep})^{-1} u_0. 
\]
Then $\theta(t)$ is written as 
\[
\theta (t) =\frac{1}{2\pi i}\int_{\Gamma} e^{\lambda t} \theta_\lambda[u_0] \, d\lambda,
\]
and $\theta_\lambda[u_0]$ satisfies 
\begin{equation}\label{5.24}
(\lambda - \Delta_D) \theta_\lambda[u_0]=\theta_0 + \mathcal{R}^{\ep}_{1,c} \bw_\lambda[u_0] \cdot \b e_2. 
\end{equation} 
Here $\Delta_D$ is the Laplace operator on $L^2(\Omega)$ under the homogeneous Dirichlet boundary condition on $\{x_2=0,1\}$, i.e., the operator on $L^2(\Omega)$ defined by 
\[
D(\Delta_D) = H^2(\Omega)\cap H_0^1(\Omega), 
\ \ \ 
\Delta_D \theta = \Delta \theta \ \ (\theta\in D(\Delta_D)). 
\]
It is well known that there exists a positive constant $\tilde\kappa_0$ such that $\{\lambda \in \mathbb{C} ; |\arg (\lambda+\tilde\kappa_0)| \le \frac{3\pi}{4} \}\subset \rho( \Delta_D )$, 
\[
\|(\lambda - \Delta_D)^{-1}h \|_2 \le \frac{C}{1+|\im\lambda|} \| h \|_2
\]
for $\lambda$ satisfying $|\arg (\lambda+\tilde\kappa_0)| \le \frac{3\pi}{4}$, and 
\[
\| e^{t\Delta_D} h\|_2^2 \le C e^{-\kappa t} \|h\|_2
\]
for $t\ge 0$, where $\kappa$ is a constant with $0 < \kappa < \tilde\kappa_0$. 

Hereafter, we fix $\kappa$ in such a way that $0 < \kappa \le \min\{\kappa_0,\tilde\kappa_0, \Lambda_1 \}$ with $\Lambda_1$ being the number given in Theorem \ref{thm4.1}. 
By \eqref{5.24}, we have 
\[
\theta_\lambda[u_0] = (\lambda - \Delta_D)^{-1} \theta_0 +\mathcal{R}^{\ep}_{1,c} (\lambda - \Delta_D)^{-1} \bw_\lambda[u_0] \cdot \b{e}_{2}, 
\]
and hence, 
\[
\theta (t)= e^{t\Delta_D} \theta_0 + \mathcal{R}^{\ep}_{1,c} \tilde\theta[u_0](t), 
\]
where 
\[
\tilde\theta[u_0](t) = \frac{1}{2\pi i}\int_{\Gamma} e^{\lambda t} (\lambda - \Delta_D)^{-1} \bw_\lambda[u_0] \,d\lambda.
\]

As for the first term $e^{t\Delta_D} \theta_0$, we have $\int_0^t e^{2\kappa s} \|e^{s\Delta_D} \theta_0\|_2^2 \,ds \le C \|\theta_0\|_2^2$. 
As for the second term $\tilde\theta[u_0](t)$, we set 
\[
\ts 
\Theta^\ep [u_0] (\eta)= \left(-\kappa + i\frac{\gamma}{\ep} - \Delta_D \right)^{-1}\bw_{-\kappa + i\frac{\eta}{\ep}}[u_0] \ \ (\gamma \in \mathbb{R}). 
\]
If $u_{0} \in R(Q^{\ep})$, we have 
\begin{equation}\label{5.25'}
\|\bw_{-\kappa +i\frac{\gamma}{\ep}}[u_0]\|_2 \le \frac{C}{1 + |\gamma| } (|||u_0|||_\ep + \ep^2\|\delx \phi_0 \|_2). 
\end{equation}
An outline of the proof of \eqref{5.25'} will be given in the end of this section. 
It then follows that  
\begin{equation}\label{5.25}
\begin{split}
\|\Theta^\ep [u_0] (\gamma) \|_2 & \le \frac{C}{1 + |\frac{\gamma}{\ep}|}\frac{1}{1 + |\gamma| } (|||u_0|||_\ep + \ep^2\|\delx \phi_0 \|_2)
\\
& \le \frac{C}{1 + |\frac{\gamma}{\ep}|}\frac{1}{1 + |\gamma| } E(u_0)^{\frac{1}{2}}
\end{split}
\end{equation}
for $\gamma \in \mathbb{R}$. 
Since $u_0\in R(Q^\ep)$, we see that $(\lambda - \Delta_D)^{-1} \bw_\lambda[u_0]$ is analytic in $\{\lambda \in \mathbb{C} ; \re\lambda \ge -\kappa_1\}$. 
Therefore, in view of \eqref{5.25}, we have $\|\Theta^\ep [u_0] (\gamma) \|_2=O(\gamma^{-2})$ as $|\gamma| \to \infty$, so we can deform $\Gamma$ to $\tilde\Gamma =\{\lambda = -\kappa+i\frac{\gamma}{\ep} ; \gamma \in \mathbb{R} \}$. 
It then follows that 
\begin{equation}\label{5.26}
\begin{split}
\tilde\theta[u_0](t)
& = \frac{1}{2\pi i}\int_{\tilde \Gamma} e^{\lambda t}(\lambda - \Delta_D)^{-1} \bw_\lambda[u_0] \, d\lambda \\
& = \frac{1}{\ep} \frac{1}{2\pi} \int_{-\infty}^{\infty} e^{-\kappa t} e^{i\frac{\gamma}{\ep}t} \Theta^\ep [u_0] (\gamma)\, d \gamma \\
& =\frac{1}{\ep}\frac{e^{-\kappa t}}{2\pi} \hat\Theta^\ep [u_0] \left({\ts -\frac{t}{\ep}} \right), 
\end{split}
\end{equation}
where $\hat\Theta^\ep [u_0](t)$ is the Fourier transform of $\Theta^\ep [u_0](\gamma)$, i.e., it is defined by $\hat\Theta^\ep [u_0](t) = \int_{-\infty}^\infty e^{-i \gamma t} \Theta^\ep [u_0](\gamma) \, d \gamma$. 
By using \eqref{5.25}, \eqref{5.26} and the Plancherel theorem, we see that 
\begin{align*}
\int_0^t e^{2\kappa s} \|\tilde{\theta}[ u_{0} ] (s)\|_2^2\,ds 
& \le \frac{1}{\ep^2} \frac{1}{(2\pi)^2} \int_{-\infty}^\infty \left\|\hat{\Theta}^\ep[u_0] \left(\frac{s}{\ep}\right) \right\|_2^2 \,ds 
\\
& = \frac{1}{\ep} \frac{1}{(2\pi)^{2} }\int_{-\infty}^\infty \| \Theta^\ep [u_0](\gamma)\|_2^2\,d\gamma 
\\
& \le \frac{C}{\ep} \int_{-\infty}^\infty \frac{E(u_0)}{\left(1+\left|\frac{\gamma}{\ep}\right|\right)^2} \,d \gamma 
\le C E(u_0).
\end{align*}
This implies that $\int_0^t e^{2\kappa s} \|\theta (s)\|_2^2\,ds \le C E(u_0)$. 
Similarly, we can show that $\int_0^t e^{2\kappa s} \|\psi (s)\|_2^2\,ds \le C E(u_0)$, and the estimate \eqref{5.23} is proved. 

Combining \eqref{5.22} and \eqref{5.23}, we have $E(u(t)) \le C e^{-2\kappa t} E[u_0]$. 
Since $u(t) = \mathscr{V}^{\ep}_{0} (t) Q^{\ep} u_{0}$, we conclude that    
\begin{equation}\label{5.27}
||| \mathscr{V}^{\ep}_{0} (t) Q^{\ep} u_{0} |||_{\ep,X^1} \le Ce^{-\kappa t} |||u_0|||_{\ep,X^1}
\end{equation}
uniformly in $t \ge 0$ and $0 < \ep \le \ep_{1}$.

We next consider the spectrum of $\mathscr{V}^{\ep}_{0} (\frac{2\pi}{a} ) Q^\ep$ on $Q^{\ep} X^1$. 
If the series $\sum_{n=0}^\infty (-1)^n (e^{-\frac{2\pi}{a}\lambda} \mathscr{V}^{\ep}_{0} (\frac{2\pi}{a} ) Q^{\ep} )^n$ absolutely converges in $\mathfrak{B}( Q^{\ep} X_{1} )$, then $[(I- e^{-\frac{2\pi}{a}\lambda}\mathscr{V}^{\ep}_{0} (\frac{2\pi}{a} ) ) Q^{\ep} ]^{-1}$ exists and it coincides with $\sum_{n=0}^\infty(-1)^n (e^{-\frac{2\pi}{a}\lambda} \mathscr{V}^{\ep}_{0} (\frac{2\pi}{a} ) Q^{\ep} )^n$. 

Let us show that $\sum_{n=0}^\infty (-1)^n (e^{-\frac{2\pi}{a}\lambda} \mathscr{V}^{\ep}_{0} (\frac{2\pi}{a} ) Q^{\ep} )^n$ absolutely converges in $\mathfrak{B}( Q^{\ep} X_{1} )$ and derive the desired uniform estimate for $[(I- e^{-\frac{2\pi}{a}\lambda}\mathscr{V}^{\ep}_{0} (\frac{2\pi}{a} ) ) Q^{\ep} ]^{-1}$. Since $( \mathscr{V}^{\ep}_{0} (\frac{2\pi}{a} ) Q^{\ep} )^n = \mathscr{V}^{\ep}_{0} (\frac{2\pi n}{a} ) Q^{\ep} $, we see from \eqref{5.27} that 
\[
\|\ts ( \mathscr{V}^{\ep}_{0} (\frac{2\pi}{a} ) Q^{\ep} )^n u_{0} \| \le Ce^{-\frac{2\pi n}{a}\kappa},
\]
where $C$ is a positive constant independent of $\ep \in (0, \ep_{1} ]$. Here $\|( \mathscr{V}^{\ep}_{0} (\frac{2\pi}{a} ) Q^{\ep} )^n\|$ denotes the operator norm of $( \mathscr{V}^{\ep}_{0} (\frac{2\pi}{a} ) Q^{\ep} )^n$ on $Q^{\ep} X^1$. 
It then follows that if $\re\lambda>-\kappa$, then  
\[
\sum_{n=0}^\infty\|(-1)^n (e^{-\frac{2\pi}{a}\lambda}{\ts \mathscr{V}^{\ep}_{0} (\frac{2\pi}{a}  ) } Q^{\ep} )^n   \| \le C\sum_{n=0}^\infty e^{-\frac{2\pi n}{a}(\re\lambda+\kappa)} = \frac{C}{1 - e^{-\frac{2\pi }{a}(\re\lambda+\kappa)}} 
\]
If $\re\lambda > -\kappa $, then $\sum_{n=0}^\infty(e^{-\frac{2\pi}{a}\lambda}  \mathscr{V}^{\ep}_{0} (\frac{2\pi}{a} ) Q^{\ep} )^n$ converges in the operator norm. We thus conclude that $( I - e^{\frac{2\pi}{a}\lambda} \mathscr{V}^{\ep}_{0} (\frac{2\pi}{a} ) ) Q^{\ep} $ has a bounded inverse on $X^1$ and $[ (I - e^{\frac{2\pi}{a}\lambda} \mathscr{V}^{\ep}_{0} (\frac{2\pi}{a} ) ) Q^{\ep} ]^{-1} = \sum_{n=0}^\infty(e^{-\frac{2\pi}{a}\lambda} \mathscr{V}^{\ep}_{0} (\frac{2\pi}{a} ) Q^{\ep})^n$ with estimate
\[
\| [ (I -e^{-\frac{2\pi}{a}\lambda} {\ts \mathscr{V}^{\ep}_{0} (\frac{2\pi}{a} )}  Q^{\ep} ]^{-1} \| \le \frac{C}{1 - e^{-\frac{2\pi}{a}(\re\lambda + \kappa)}}.
\]
Setting $\kappa_1=\kappa$, we obtain the desired result. This completes the proof. 
\hfill$\square$

\vspace{1ex}
To prove Theorem \ref{thm5.6}, we next consider the following equation 
\begin{equation}\label{5.28} 
(I - e^{ - \frac{2\pi}{a^{\ep} } \lambda } \mathscr{V}^{\ep}_{0} ({\ts \frac{2\pi}{a} })  ) u = F, 
\end{equation}
where $F\in X^1$ is a given function. 

\vspace{1ex}
\begin{lem}\label{lem5.7} 
{\rm (i)} If $e^{\frac{2\pi}{ a^{\ep} }\lambda} \neq 1$ and $\re{\lambda} > - \frac{a^{\ep}}{a}\kappa_1$, then for any $F\in X^1$ there exists a unique solution $u \in X^1$ of \eqref{5.28}. 
The solution $u$ is given by 
\[
u=\frac{1}{1- e^{-\frac{2\pi}{a^{\ep}}\lambda}}P^\ep F + \left[ (I - e^{-\frac{2\pi}{a^{\ep}} \lambda} \mathscr{V}^{\ep}_{0}({\ts \frac{2\pi}{a} } ) )Q^\ep \right]^{-1} Q^{\ep} F, 
\]
and satisfies the estimate 
\[
||| u |||_{\ep, X^1} \le \frac{1}{{ | 1- e^{-\frac{2\pi}{a^{\ep} }\lambda} | }} ||| P^\ep F |||_{\ep,X^1} + \frac{C}{1-e^{-\frac{2\pi}{a^{\ep} }(\re\lambda + \frac{a^{\ep} }{a} \kappa_1)} }||| Q^\ep F|||_{\ep, X^1}. 
\] 

\vspace{1ex}
{\rm (ii)} Let $e^{\frac{2\pi}{a^{\ep} }\lambda}=1$. Then the equation \eqref{5.28} is solvable if and only if $P^\ep F=0$, i.e., $(F, u^{\ep *}_{+} )_{\ep} = (F, u^{\ep *}_{-} )_{\ep}=0$. If this condition is satisfied, then $F=Q^\ep F$ and the equation \eqref{5.28} is uniquely solvable in $Q^\ep X^1$; $u$ is given by 
\[
u=\left[ (I - e^{-\frac{2\pi}{a^{\ep}} \lambda } \mathscr{V}^{\ep}_{0} ({\ts \frac{2\pi}{a} } ) )
Q^{\ep} \right]^{-1} F,
\]
and satisfies the estimate 
\[
||| u |||_{\ep, X^1} \le \frac{C}{1-e^{-\frac{2\pi}{a^{\ep} }(\re\lambda + \frac{ a^{\ep} }{a} \kappa_1)} }||| F |||_{\ep, X^1}. 
\] 
\end{lem}

\vspace{1ex}
\noindent
{\bf Proof.} 
We first recall that 
\begin{equation}\label{5.29}
P_\pm^\ep \mathscr{V }^{\ep}_{\omega} (t) u = \mathscr{V}^{\ep}_{\omega}(t) P_\pm^\ep u = e^{\pm \frac{iat}{1 + \omega} } P_\pm^\ep u.
\end{equation}

Applying $P^\ep$ and $Q^\ep$ to \eqref{5.28}, we obtain
\begin{align}\label{5.30}
\begin{cases}
(1 - e^{-\frac{2\pi}{a^{\ep} }\lambda} ) P^\ep u = P^\ep F,\\
( I - e^{-\frac{2\pi}{a^{\ep} } \lambda }  \mathscr{V}^{\ep}_{0}({\ts \frac{2\pi}{a} } ) ) Q^\ep u = Q^\ep F.
\end{cases}
\end{align}
If $ e^{\frac{2\pi}{a^{\ep} }\lambda} \neq 1$, then
\[
P^\ep u = \frac{1}{ 1- e^{-\frac{2\pi}{a^{\ep} }\lambda}} P^\ep F.
\]
On the other hand, by Lemma \ref{lem5.5}, if $\re{\lambda} > - \frac{ a^{\ep} }{a}\kappa_{1}$, then 
\begin{equation}\label{5.31}
Q^\ep u = \left[ ( I- e^{-\frac{2\pi}{a^{\ep} }\lambda} \mathscr{V}^{\ep}_{0}({\ts \frac{2\pi}{a} } ) ) Q^\ep\right]^{-1} Q^{\ep} F, 
\end{equation} 
and 
\[
||| Q^\ep u|||_{\ep,X^1} \le \frac{C}{1-e^{ -\frac{2\pi}{a^{\ep} }(\re\lambda+ \frac{a^{\ep} }{a} \kappa_1)}}||| Q^\ep F|||_{1,X^1}.
\]
We thus prove (i). 

If $ e^{\frac{2\pi}{a^{\ep} }\lambda} = 1$, the first equation of \eqref{5.30} requires $P^\ep F =0$, i.e., $(F, u^{\ep *}_{+} )_{\ep} = (F, u^{\ep *}_{-} )_{\ep}=0$; and in this case $F = Q^\ep F$. Therefore, if $\re{\lambda} > \frac{a^{\ep} }{a}\kappa_{1}$, then \eqref{5.28} is solvable and a solution $u$ is given by \eqref{5.31} and it is unique under the condition $P^\ep u =0$, i.e., $u \in Q^\ep X^1$. 
This completes the proof. \hfill$\square$

\vspace{1ex}
We now prove Theorem \ref{thm5.6}. 

\vspace{1ex}
\noindent
{\bf Proof of Theorem \ref{thm5.6}.} 
The assertion (i) follows from (ii). In fact, since $L^{\ep}_{\mathcal{R}_{1,c}^\ep} u_\pm^\ep =\mp ia^{\ep} u_\pm^\ep$, we see that $B^{\ep} (e^{i(k \pm 1) a t} u_\pm^\ep ) = i (k \pm 1)a^{\ep} u_\pm^\ep \mp ia^{\ep} u_\pm^\ep = ika^{\ep} u_\pm^\ep$. Therefore, $ika^{\ep} $ is an eigenvalue for each $k\in \mathbb{Z}$. From the representation formula of $(\lambda + B^{\ep})^{-1}$ we see that singularities of $(\lambda + B^{\ep})^{-1}$ in $\Sigma_1$ are given by the points $\lambda$ satisfying $e^{\frac{2\pi}{a^{\ep} }\lambda} = 1$; so these points take the forms $\lambda=ika^{\ep} $ $(k\in \mathbb{Z})$. At each $ika^{\ep} $, $(1- e^{\frac{2\pi}{a^{\ep} }\lambda})^{-1}=O((\lambda - ika^{\ep} )^{-1})$, which implies that $(\lambda + B^{\ep})^{-1} = O((\lambda - ika^{\ep} )^{-1})$ as $ \lambda \to ika^{\ep}  $. This shows that each $ika^{\ep} $ is a semisimple eigenvalue of $B^{\ep}$.  
Therefore, (i) follows from (ii).

Let us prove (ii). For a given $F\in \mathcal{Y}_1$, we consider the system
\begin{equation}\label{5.32}
(\lambda + B^{\ep})u = F, \, \ \ u \in \mathcal{Y}_a. 
\end{equation}
If $u \in \mathcal{Y}_a$, then $u(0)\in X^1$. It then follows from Lemma \ref{lem5.1} with $\beta = \frac{a^{\ep} }{a}$ that solution of \eqref{5.32} is written in the form 
\begin{equation}\label{5.33}
u(t) =  e^{-\frac{a}{a^{\ep} }\lambda t} \mathscr{V}^{\ep}_{0}(t) u(0) + \frac{a}{a^{\ep} } \int_0^t e^{-\frac{a}{a^{\ep} } \lambda (t-s) }  \mathscr{V}^{\ep}_{0} ( t - s ) F(s) \,ds. 
\end{equation}
The periodicity condition $u\left(\frac{2\pi}{a}\right) = u(0) \in X^1$ requires
\begin{equation}\label{5.34}
(I- e^{-\frac{2\pi}{a^{\ep}} \lambda} \mathscr{V}^{\ep}_{0}( {\ts \frac{2\pi}{a} } ) ) u(0) = \frac{a}{a^{\ep} } \int_0^{\frac{2\pi}{a}} e^{-\frac{a}{a^{\ep} }\lambda (\frac{2\pi}{a}-s )} \mathscr{V}^{\ep}_{0}({\ts \frac{2\pi}{a} - s } ) F(s)\,ds.  
\end{equation}
If one could find $u(0)\in X^1$ satisfying \eqref{5.34}, then the function $u(t)$ in \eqref{5.33} with this $u(0)$ would give a solution of \eqref{5.32}. So we look for $u(0)\in X^1$ satisfying \eqref{5.34}. 

Let $\tilde{F} = \frac{a}{a^{\ep} } \int_0^{\frac{2\pi}{a}} e^{-\frac{a}{a^{\ep} }\lambda (\frac{2\pi}{a}-s )} \mathscr{V}^{\ep}_{0}({\ts \frac{2\pi}{a} - s } ) F(s)\,ds$. By Lemma \ref{lem5.1} with $u_0 = 0$ we see that $\tilde{F}\in X^1$. 
Furthermore, since $e^{\mp ias} P_\pm^\ep F = (F, z_\pm^{\ep*}(s))_\ep u_\pm$,  we see from \eqref{5.29} that  
\begin{align*}
P_+^\ep \tilde{F} 
              &= \frac{a}{a^{\ep} } \int_0^{\frac{2\pi}{a}} e^{-\frac{a}{a^{\ep} }\lambda (\frac{2\pi}{a}-s )} P^{\ep}_{+} \mathscr{V}^{\ep}_{0}({\ts \frac{2\pi}{a} - s } ) F(s)\,ds.  \\
              &= \frac{a}{a^{\ep}} \int_0^{\frac{2\pi}{a}} e^{-\frac{a}{a^{\ep} }\lambda (\frac{2\pi}{a}-s)+ia(\frac{2\pi}{a}-s) } P_+^\ep F(s) \,ds \\
              &= \frac{a}{a^{\ep} }\int_0^{\frac{2\pi}{a}} e^{-\frac{a}{a^{\ep} } \lambda (\frac{2\pi}{a}-s) } (F(s), z_+^{\ep*}(s) )_\ep u^{\ep}_{+} \,ds \\ 
              & = \frac{2\pi}{a^{\ep} }\left[ e^{-\frac{a}{a^{\ep}} \lambda (\frac{2\pi}{a}-s ) } F\right]_{+,\ep} u^{\ep}_{+}.
\end{align*}
Similarly, we obtain $P_-^\ep \tilde{F} = \frac{2\pi}{a^{\ep} } \left[ e^{-\frac{a}{a^{\ep} } \lambda (\frac{2\pi}{a}-s ) }F\right]_{-,\ep} u_-^\ep$, and hence, 
\begin{equation}\label{5.35}
P^\ep \tilde{F} = \frac{2\pi}{a^{\ep} }\left[ e^{-\frac{a}{a^{\ep} }\lambda (\frac{2\pi}{a}-s) }F\right]_{+,\ep}u_{+}^{\ep} + \frac{2\pi}{a^{\ep} }\left[ e^{-\frac{a}{a^{\ep}} \lambda (\frac{2\pi}{a}-s) } F\right]_{-,\ep}u_{-}^{\ep}.
\end{equation}
It then follows from Lemma \ref{lem5.7} that if $ e^{\frac{2\pi}{a^{\ep} }\lambda}\neq1$ and $\re{\lambda} > - \frac{a^{\ep}}{a}\kappa_{1}$, then
\[
P^\ep u(0) = \frac{ 2\pi}{a^{\ep} (1- e^{-\frac{2\pi}{a^{\ep} }\lambda})} \left(\left[ e^{-\frac{a}{a^{\ep} }\lambda (\frac{2\pi}{a}-s) } F \right]_{+,\ep}u_{+}^{\ep} + \left[ e^{-\frac{a}{a^{\ep} }\lambda (\frac{2\pi}{a}-s) } F \right]_{-,\ep}u_{-}^{\ep}\right)
\]
and 
\[
Q^\ep u(0) = \frac{a}{ a^{\ep} } \left[(I - e^{ -\frac{2\pi}{a^{\ep} } \lambda} \mathscr{V}^{\ep}_{0}({\ts \frac{a}{a^{\ep} } }) ) Q^{\ep} \right]^{-1}Q^\ep \int_0^{\frac{2\pi}{a}} e^{- \frac{a}{a^{\ep} } \lambda (\frac{2\pi}{a}-s)} \mathscr{V}^{\ep}_{0}( {\ts \frac{2\pi}{a} - s } )F(s)\,ds.
\]
We thus find that $u(0) = P^\ep u(0) + Q^\ep u(0)$ satisfies \eqref{5.34}, and conclude that if $ e^{\frac{2\pi}{a^{\ep} }\lambda}\neq1$ and $\re{\lambda} > - \frac{a^{\ep} }{a} \kappa_{1}$, the solution $u(t)$ of \eqref{5.32} is given by the formula described in (ii). 
This proves (ii). 

The assertion (iii) is a direct consequence of the semi-simplicity of the eigenvalue $0$ which is assured in (i). 

As for (iv), assume that the problem $B^{\ep}u = F$ with $[u]_{+,\ep} = [u]_{-,\ep} = 0$ has a unique solution. Since $[B^{\ep} u]_{+,\ep} = [B^{\ep} u]_{-,\ep} =0$, we have $[F]_{+,\ep} = [F]_{-,\ep} = 0$.

Conversely, assume that $[F]_{+,\ep} = [F]_{-,\ep} = 0$. We first observe that \eqref{5.35} holds for $\lambda = 0$ with $\tilde{F} = \frac{a}{a^{\ep}} \int_0^{\frac{2\pi}{a}} \mathscr{V}^{\ep}_{0}({\ts \frac{2\pi}{a}-s } ) F(s)\,ds$. It then follows that 
\[
P^\ep\left[\frac{a}{a^{\ep}} \int_0^{\frac{2\pi}{a}} \mathscr{V}^{\ep}_{0}({\ts \frac{2\pi}{a}-s } ) F(s)\,ds\right] = \frac{2\pi}{a^{\ep} } \left([F]_{+,\ep}u_{+}^{\ep} + [F]_{-,\ep}u_{-}^{\ep}\right) = 0, 
\]
and hence,
\[
Q^\ep \int_0^{\frac{2\pi}{a}} \mathscr{V}^{\ep}_{0}({\ts \frac{2\pi}{a}-s } ) F(s) \,ds  =  \int_0^{\frac{2\pi}{a}} \mathscr{V}^{\ep}_{0}({\ts \frac{2\pi}{a}-s } ) F(s) \,ds.
\]
Therefore, we see that $u(0)$ is given by
\begin{equation}\label{5.36}
u(0) = \frac{a}{a^{\ep}} \left[ (I -  \mathscr{V}^{\ep}_{0}({\ts \frac{2\pi}{a} } ) )Q^\ep \right]^{-1}\int_0^{\frac{2\pi}{a}} \mathscr{V}^{\ep}_{0}({\ts \frac{2\pi}{a}-s } ) F(s)\,ds
\end{equation}
satisfies \eqref{5.34} with $\lambda = 0$. We thus obtain a solution $\tilde u(t)$ by taking $u(0)$ given by \eqref{5.36} in the formula \eqref{5.33} with $\lambda = 0$: 
\begin{align*}
\tilde u(t) & =  \frac{a}{a^{\ep}} \mathscr{V}^{\ep}_{0}(t) \left[ (I -  \mathscr{V}^{\ep}_{0}({\ts \frac{2\pi}{a} } ) )Q^\ep \right]^{-1}\int_0^{\frac{2\pi}{a}} \mathscr{V}^{\ep}_{0}({\ts \frac{2\pi}{a}-s } ) F(s)\,ds\\ 
& \quad + \frac{a}{a^{\ep}  }\int_0^t  \mathscr{V}^{\ep}_{0}(t - s) F(s)\,ds F(s)\, ds.
\end{align*}
By using \eqref{5.29} and the fact $\mathcal{P}_\pm^\ep G = e^{\pm iat} \left(\frac{a}{2\pi} \int_0^{\frac{2\pi}{a} } e^{\mp ias} P_\pm^\ep G(s)\,ds \right)$, we have 
\[
\mathcal{P}_{\pm}^\ep \left ( \int_0^t\mathscr{V}^{\ep}_{0}(t - s)  F(s)\, ds \right) = \frac{2\pi}{a} \mathcal{P}_{\pm}^\ep(F) - \mathcal{P}_{\pm}^\ep (sF(s)) = - \mathcal{P}_{\pm}^\ep (sF(s)) 
\]
if $[F]_{+,\ep} = [F]_{-,\ep} = 0$. Setting $u = \tilde u - \frac{a}{a^{\ep} } \mathcal{P}^\ep (sF(s))$, we see that $u$ is a solution of $B^{\ep} u = F$ and $u \in \mathcal{Q}^\ep \mathcal{X}_a$ ; and $u$ is written in the formula given in (iv).  
As for the uniqueness, one can show that $B^{\ep}u = 0$ and $[u]_{\pm,\ep} = 0$ imply $u=0$. In fact, let $B^{\ep}u = 0$. Since $N(B^{\ep}) = {\rm span}\,\{z_+^\ep, z_-^\ep \}$, we have $u=c_+ z_+^\ep + c_- z_-^\ep$. But since $[u]_{\pm,\ep} = 0$, we see that $c_+ = c_- = 0$, and hence, $u =0$. Using the formula in (iii), together with Lemmas \ref{5.1} and \ref{5.5}, we obtain the estimate $
||| u |||_{\ep,\mathcal{Y}_a} \le C||| F |||_{\ep,\mathcal{X}_a}$. This completes the proof. \hfill$\square$

\vspace{2.5ex}
We next prove Lemma \ref{lem5.8}. 

\vspace{1ex} 
\noindent
{\bf Proof of Lemma \ref{lem5.8}.} We first observe that $\mathcal{P}_\pm^\ep \delt \subset \delt \mathcal{P}_\pm^\ep = \pm ia \mathcal{P}_\pm^\ep$. Therefore, $\mathcal{P}^\ep \mathcal{X}_a$ is invariant under $\delt$ and so is $\mathcal{Q}^\ep \mathcal{Y}_a$. 
We thus deduce that $\mathcal{Q}^\ep \mathcal{Y}_a$ is invariant under $B^{\ep}( \omega)$. 

For a given $F \in \mathcal{X}_{a}$, we consider the equation  
\begin{equation}\label{5.37}
B^{\ep}(\omega) u = F.  
\end{equation}  

If $\omega$ is small enough, say, $|\omega| \le \frac{1}{4}$, then $\sigma( -\frac{a}{a^{\ep} (1 + \omega)}L^{\ep}_{ \mathcal{R}^{\ep}_{1,c} } ) \cap \{\lambda\in \mathbb{C}; \re\lambda > -\frac{a}{a^{\ep}(1+\omega) } \Lambda_{1} \} = \{ \frac{ia}{1 + \omega}, -\frac{ia}{1 + \omega} \}$. 
As in Theorem \ref{thm5.6}, this implies that $\sigma(- B^{\ep}(\omega)) \cap \{\lambda\in \mathbb{C}; \re\lambda > -\frac{a^{\ep} }{a} \kappa_{1} \} = \cup_{k \in \mathbb{Z} } \left\{ i a^{\ep} ( (1 + \omega ) k + 1 ),  i a^{\ep} ( (1 + \omega ) k - 1 )  \right\}$. 
Therefore, if $|\omega| \le \frac{1}{4}$ and $\omega \neq 0$, then $0 \in \rho(- B^{\ep}( \omega))$, 
and hence, for any $F\in \mathcal{X}_a$, the equation \eqref{5.37} has a unique solution $u\in \mathcal{Y}_a$. 

By \eqref{5.29}, we have $P_\pm^{\ep}(I - \mathscr{V}^{\ep}_{\omega}( {\ts \frac{2\pi}{a} }) )=(1-e^{\pm \frac{2\pi i}{1+\omega}}) P^{\ep}_{\pm}$. 
Furthermore, by Lemma \ref{lem5.5}, if $|\omega| \le \frac{1}{4}$, $||| \mathscr{V}^{\ep}_{\omega}(t) Q^\ep u_0 |||_{\ep, X^1} \le C e^{- \frac{\kappa_1 }{1+ \omega}t } ||| u_0 |||_{\ep, X^1} $. 
Therefore, as in the proof of Theorem \ref{thm5.6}, if $|\omega| \le \frac{1}{4}$ and $\omega \neq 0$, then $I -  \mathscr{V}^{\ep}_{\omega}({\ts \frac{2\pi}{a} } )  $ has a bounded inverse on $X^1$, and the solution $u$ of \eqref{5.37} is represented as 
\begin{equation}\label{5.38}
u (t)  =  u_1(t) + u_2(t) + u_3(t), 
\end{equation} 
where 
\begin{align*}
u_1(t) & =  \frac{2\pi}{a^{\ep} ( 1 + \omega )} \sum_{ j= +, -} \frac{  e^{ j \frac{2\pi i}{1+\omega} } }{ 1 - e^{ j \frac{2\pi i}{1+\omega} } } \left[ e^{ j \frac{ia \omega }{1 + \omega} s } F \right]_{j, \ep} e^{ j \frac{ia}{1 + \omega} t } u_j^\ep, \\
u_2(t) & = \frac{a}{a^{\ep} (1 + \omega)} \mathscr{V}^{\ep}_{\omega} (t) \left[ (I -  \mathscr{V}^{\ep}_{\omega}({\ts \frac{2\pi}{a} } ) ) Q^\ep \right]^{-1}
Q^\ep \int_0^{\frac{2\pi}{a}} \mathscr{V}^{\ep}_{\omega} ({\ts \frac{2\pi}{a} -s } ) F(s)\,ds, \\
u_3(t) & = \frac{a}{a^{\ep}(1+\omega)} \int_0^t  \mathscr{V}^{\ep}_{\omega}( t - s ) F(s)\,ds.
\end{align*}
As observed at the beginning of the proof, the solution $u \in \mathcal{Q}^\ep \mathcal{Y}_a$ if $F \in \mathcal{Q}^\ep \mathcal{X}_a$. 

We next establish a uniform estimate of $u$ with respect to $0< \ep \le \ep_{1}$ and $|\omega| \le \frac{1}{4}$. To this end, let us compute the $ \mathcal{Q}^\ep $ part of the right-hand side of \eqref{5.38} when $F \in \mathcal{Q}^\ep \mathcal{X}_a$. 
By using \eqref{5.29} and the fact $\mathcal{P}_\pm^\ep G = e^{\pm iat} \left(\frac{a}{2\pi} \int_0^{\frac{2\pi}{a} } e^{\mp ias} P_\pm^\ep G(s)\,ds \right)$, we find that 
\begin{align*}
\mathcal{P}_\pm^\ep u_1 & = \mp \frac{1}{ia^{\ep} \omega} (e^{\mp \frac{2\pi i }{1 + \omega}}-1) ^{-1} ( e^{\mp \frac{2\pi i \omega}{1 + \omega}}-1) \mathcal{P}_\pm^\ep \left(e^{ \pm \frac{ia \omega }{1+\omega}s } F \right), \\ 
\mathcal{P}_\pm^\ep u_2 & = 0, \\ 
\mathcal{P}_\pm^\ep u_3 & = \mp \frac{1}{ia^{\ep} \omega} \mathcal{P}_\pm^\ep \left( \left( e^{ \mp \frac{2\pi i \omega}{1 + \omega} \pm \frac{ia \omega }{1+\omega}s } - 1 \right) F \right). 
\end{align*} 
In particular, we obtain $\mathcal{P}^\ep u = \mathcal{P}^\ep (u_1 + u_2 + u_3) = \frac{1}{ia^{\ep} \omega} \left(\mathcal{P}_+^\ep F - \mathcal{P}_-^\ep F \right)$, so $\mathcal{P}^\ep (u_1 + u_2 + u_3) = 0 $ when $F \in \mathcal{Q}^\ep \mathcal{X}_a$. 
Since $\mathcal{Q}^\ep = I - \mathcal{P}^\ep$, we find that if $F \in \mathcal{Q}^\ep \mathcal{X}_a$, then 
\[
u = \mathcal{Q}^\ep u  = u_1 + u_2 + u_3.
\]
From the proofs of Lemmas \ref{lem5.1} and \ref{lem5.5}, it follows that $||| u_j |||_{\mathcal{Y}_a} \le C ||| F |||_{\mathcal{X}_a}$ for $j=2,3$ uniformly for $0< \ep \le \ep_{1}$ and $|\omega| \le \frac{1}{4}$. Furthermore, since $F \in \mathcal{Q}^\ep \mathcal{X}_a$, we have $[F]_{\pm,\ep} = 0$, and therefore, 
\[
\frac{e^{ \pm \frac{2\pi i}{1+\omega} } }{1 - e^{ \pm \frac{2\pi i}{1+\omega} } } \left[ e^{ \pm \frac{ia \omega }{1 + \omega} s } F \right]_{\pm , \ep} 
= \left[ \frac{e^{ \pm \frac{2\pi i}{1+\omega} } }{1 - e^{ \pm \frac{2\pi i}{1+\omega} } } \left( e^{ \pm \frac{ia \omega }{1 + \omega} s } - 1 \right) F \right]_{\pm , \ep} .
\]
Noting that $e^{ \pm \frac{2\pi i}{1+\omega} } = e^{ \pm \frac{2\pi i}{1+\omega} \mp 2\pi i} = e^{ \mp \frac{2\pi \omega i}{1+\omega} }$ and  
\[
\frac{1}{1 - e^{ \pm \frac{2\pi i}{1+\omega} }} \left( e^{ \pm \frac{ia \omega }{1 + \omega} s } - 1 \right) 
= \frac{ \pm \frac{ia \omega }{1 + \omega} s }{ 1 - e^{ \mp \frac{2\pi \omega i}{1+\omega} } } 
\int_0^1 e^{ \pm \frac{ia \omega }{1 + \omega} s \sigma} \,d\sigma \to \frac{as}{2\pi
} 
\]
as $\omega\rightarrow0$ uniformly in $0 \le s \le \frac{2\pi}{a}$, we see that $||| u_1 |||_{\mathcal{Y}_a} \le C ||| F |||_{\mathcal{X}_a}$ for $0< \ep \le \ep_{1}$ and $0< |\omega| \le \frac{1}{4}$. 
This estimate also holds for $\omega =0$. This completes the proof. 
\hfill$\square$.  
 
\begin{rem}\label{rem5.9} 
By \eqref{5.29}, we have $P_\pm^{\ep}(\mu - \mathscr{V}^{\ep}_{\omega}( {\ts \frac{2\pi}{a} }) )=(\mu - e^{\pm \frac{2\pi i}{1+\omega}}) P^{\ep}_{\pm}$ for any $\mu \in \mathbb{C}$. 
Furthermore, by Lemma {\rm \ref{lem5.5}}, if $|\omega| \le \frac{1}{4}$, $||| \mathscr{V}^{\ep}_{\omega}(t) Q^\ep u_0 |||_{\ep, X^1} \le C e^{- \frac{\kappa_1 }{1+ \omega}t } ||| u_0 |||_{\ep, X^1} $. 
Therefore, for any $r > 0$, there exists a positive constant $\omega_{r} = O(r)$ as $r \to 0$ such that if $|\omega| \le \omega_{r}$ then $\mu -  \mathscr{V}^{\ep}_{\omega}({\ts \frac{2\pi}{a} } ) $ has a bounded inverse on $X^1$ for $\mu \in \mathbb{C}$ satisfying $|\mu - 1| \ge r$ and $|\mu| \ge e^{-\frac{3}{4}\kappa_{1}\frac{2\pi}{a} }$ and $(\mu -  \mathscr{V}^{\ep}_{\omega}({\ts \frac{2\pi}{a} }) )^{-1} $ satisfies the estimate 
\[
||| (\mu -  \mathscr{V}^{\ep}_{\omega}({\ts \frac{2\pi}{a} } ) )^{-1} |||_{\ep, X^{1} } \le C\left(\frac{1}{r} + \frac{1}{|\mu|} \right) ||| F |||_{\ep, X^{1}} 
\]
uniformly for $0 < \ep \le \ep_{1}$. 

%Similarly to above, one can also obtain the following solution formula for $( \lambda + B^{\ep}(\omega) ) u = F$ for $\lambda$ satisfying $\re\lambda \ge - \frac{a^{\ep} }{a} \kappa_{1}$ and $\lambda \neq i a^{\ep} ((1 + \omega ) k \pm 1)$ $(k \in \mathbb{Z})$: 
%\begin{equation}\label{5.39}
%u (t)  =  u_1(t) + u_2(t) + u_3(t), 
%\end{equation} 
%where 
%\begin{align*}
%u_1(t) & =  \frac{2\pi}{a^{\ep} ( 1 + \omega )} \sum_{ j= +, -} \frac{1}{ e^{ \frac{2\pi (\lambda - j i a^{\ep} ) }{ a^{\ep} ( 1 + \omega )} } - 1 } \left[ e^{ \frac{ a( \lambda + j i a^{\ep} \omega ) }{ a^{\ep} (1 + \omega) } s } F \right]_{j, \ep} e^{ - \frac{ a ( \lambda  - j i a^{\ep} \omega ) }{ a^{\ep} (1 + \omega) } t } u_j^\ep, \\
%u_2(t) & = \frac{a}{a^{\ep} (1 + \omega)} \mathscr{V}^{\ep}_{\omega} (t) \left[ \left( I -  e^{ -\frac{2 \pi \lambda }{ a^{\ep} ( 1 + \omega ) } } \mathscr{V}^{\ep}_{\omega}({\ts \frac{2\pi}{a} } ) \right) Q^\ep \right]^{-1} \\[2ex]
%& \quad \quad \cdot Q^\ep \int_0^{\frac{2\pi}{a}} e^{ -\frac{ a \lambda }{ a^{\ep} ( 1 + \omega ) } \left( \frac{ 2 \pi }{a} - s \right) } \mathscr{V}^{\ep}_{\omega} ({\ts \frac{2\pi}{a} -s } ) F(s)\,ds, \\
%u_3(t) & = \frac{a}{a^{\ep}(1+\omega)} \int_0^t  e^{ -\frac{ a \lambda }{ a^{\ep} ( 1 + \omega ) } ( t - s ) } \mathscr{V}^{\ep}_{\omega}( t - s ) F(s)\,ds.
%\end{align*}
%This representation will be employed to study the stability of the bifurcating time periodic solutions. 
\end{rem}

In the end of this section we give an outline of the proof of \eqref{5.25'} 

\vspace{2ex}
\noindent
{\bf Outline of Proof of (\ref{5.25'}).} 
As was observed in \cite{Kagei-Nishida, Kagei-Nishida-Teramoto, Kagei-Teramoto}, the spectrum of $- L^{\ep}_{\mathcal{R}^{\ep}_{1,c}}$ is decomposed into two parts; one lies in a region with $|\im \lambda| = O(\ep^{-1})$ and the other one lies in a region $|\im\lambda|=O(1)$ as $\ep \to 0$; the latter part consists of eigenvalues which are given by perturbations of eigenvalues of $- \mathbb{L}_{\mathcal{R}_{1,c}}$. Based on this observation, we estimate the $Q^{\ep}$ part of the resolvent $(\lambda + L^{\ep}_{\mathcal{R}^{\ep}_{1,c}})^{-1}$ by dividing the region of $\lambda$ into two parts to obtain the inequality \eqref{5.25'}. 

\vspace{1ex}
We first consider the part $|\lambda| \le C \ep^{-1}$ where $(\lambda + L^{\ep}_{\mathcal{R}^{\ep}_{1,c}})^{-1}$ can be regarded as a {\it perturbation} of $(\lambda + \mathbb{L}_{\mathcal{R}_{1,c}})^{-1}$ if $\ep$ is sufficiently small. We set 
\[
||| u |||_{\ep, X} = \sqrt{\ep^{2} \|\phi\|_{2}^{2} + \ep^{4} \|\nabla \phi \|_{2}^{2} + \|\bu\|_{2}^{2} } \quad \quad (u = \trans(\phi, \bu)).
\]

\begin{lem}\label{lemA.3'} 
There exist positive constants $C_{0}$ and $C$ such that if $\lambda \in \Sigma\cap \{\lambda; |\lambda| \le C_{0}\ep^{-1}\}$ and $0< \ep \le \ep_{1}$, then the problem $(\lambda + L^{\ep}_{\mathcal{R}^{\ep}_{1,c}} ) u =Q^{\ep}F$ has a unique solution $u_{\lambda} = \trans(\phi_{\lambda}, \bu_{\lambda})$ in $D(L^{\ep}_{\mathcal{R}^{\ep}_{1,c}}) \cap R(Q^{\ep})$ and $u_{\lambda} = \trans(\phi_{\lambda}, \bu_{\lambda})$ satisfies the estimates  
\begin{align*}
\ep \|\phi_{\lambda} \|_{2} + \ep^{2} \|\nabla \phi_{\lambda} \|_{2} & \le C ||| F |||_{\ep, X},  \\[1ex]
\|\bu_{\lambda} \|_{2} & \le \frac{C}{ 1 + |\lambda| } ||| F |||_{\ep, X}. 
\end{align*}
\end{lem} 

Lemma \ref{lemA.3'} can be proved in a similar perturbation argument to that given in \cite[Section 4.1]{Kagei-Teramoto} by using the norm $||| u |||_{\ep, X}$. 

\vspace{1ex}
We next consider the part $|\im \lambda|  \ge C_{0}\ep^{-1}$. We recall that the Poincar\'e inequality $\|\nabla \bu\|_{2}\ge c_{P} \|\bu\|_{2}$ holds for $\bu \in \b X^{1}$ with some constant $c_{P}$ depending only on $\alpha$.  

\begin{lem}\label{lemA.3}
Let $\kappa_2$ and $\gamma_*$ be given positive numbers. 
There exist a positive constant $\tilde{\kappa}_{1}$ depending only on $c_{P}^2\Pr$ such that if $0<\ep \le \ep_{1}$, then  
\[
\rho(-L^{\ep}_{\mathcal{R}_{1,c}^\ep})
\supset 
\left\{
\lambda=\kappa+i\frac{\gamma}{\ep}; -\tilde{\kappa}_{1} \le \kappa\le\kappa_2, \, |\gamma| \ge \gamma_{*} \right\}.
\] 
Furthermore, if $\trans(\phi_{\lambda}, \bw_{\lambda}, \theta_{\lambda}, \psi_{\lambda}) = (\lambda+L^{\ep}_{\mathcal{R}_{1,c}^\ep})^{-1}F$ with $F=\trans(f,\b F)$ and $\lambda = \kappa + i \frac{\gamma}{\ep}$, then the estimate  
\[
(\gamma+ c_{P})\|\bw_{\lambda}\|_2+\|\nabla\bw_{\lambda}\|_2\le C ||| F |||_{\ep, X} 
\]
holds for uniformly for $0<\ep \le \ep_{1}$ and $\lambda=\kappa+i\frac{\gamma}{\ep}$ with $-\tilde{\kappa}_{1} \le \kappa \le \kappa_2$ and $\gamma \ge \gamma_{*}$. 
%\[
%\ep\|\delx^{2}\bu\|_{2}+\ep\|\delx \phi\|_{2}
%\le 
%C\left\{\|\b F\|_{2}+\ep\|f\|_{H^{1}}\right\}.
%\]
\end{lem} 

Lemma \ref{lemA.3} is proved in a similar manner to the proofs of \cite[Proposition 6.5]{Kagei-Nishida} and \cite[Proposition 3.5]{Kagei-Nishida-Teramoto}. In fact, by a similar argument as those in the proofs there, one can obtain the estimate 
\begin{align*}
& \frac{1}{8  \Pr } ( 16\kappa\gamma + \gamma^{3} -\ep^2 \gamma \kappa ) \| \b{w} \|_{2}^{2}  + \frac{1}{2} (\gamma - 2 \ep^{2} \gamma |\kappa| ) \| \nabla \b{w} \|_{2}^{2}\\
& \quad \le C ( \ep^{2} | \kappa | + \ep^{2} |\kappa| \gamma^{-1} + \ep ) \| \b{w} \|_{2}^{2} \\
& \quad \quad + C \frac{ \ep }{ \gamma } (\ep \gamma + 1)( \ep | \kappa | + \gamma) ( \|h\|_{2} + \| k \|_{2} ) \|\b{w}\|_{2} \\
& \quad \quad + \ep \im( G_{\lambda}, \b{w} )_{L^2} + \ep^{2} \gamma \re ( G_{\lambda}, \b{w} )_{L^2},
\end{align*}
where $G_{\lambda} = \lambda \b{g} - \nabla f$ and $C$ is a positive constant depending only on $\eta_{2}$, $\mathcal{R}_{1,c}$, $\mathcal{R}_{2,*}$, $\mathcal{R}_2^*$, and  $\tau$. 
Noting that $\im( G_{\lambda}, \b{w} )_{L^2} = \im\{\lambda(\b{g}, \bw)_{L^2} + (f, \div\b{w})_{L^2} \}$, one can obtain the desired estimate for $\b{w}_{\lambda}$ in Lemma \ref{lemA.3} similarly to the arguments given in \cite{Kagei-Nishida, Kagei-Nishida-Teramoto}. 

\vspace{1ex}
The estimate \eqref{5.25'} follows from Lemmas \ref{lemA.3'} and \ref{lemA.3} by a suitable choice of $\gamma_{*} >0$. 
\hfill$\square$   

\vspace{2ex}
\noindent
{\bf Acknowledgements.}
Y. Kagei was partly supported by JSPS KAKENHI
Grant Numbers 16H03947, 16H06339 and 20H00118. T. Nishida is partly supported by JSPS KAKENHI Grant Number 20K03699.

\end{document}